\documentclass[11pt]{article}

\usepackage[T1]{fontenc}
\usepackage[utf8]{inputenc}
\usepackage{amsmath,amssymb,epsfig,bbm}
\usepackage{stmaryrd}
\usepackage[colorlinks]{hyperref}
\usepackage{comment}
\usepackage{color}
\usepackage{textcomp}

\usepackage[acronym,toc]{glossaries}
\usepackage{glossary-longragged}
\usepackage{glossary-superragged}
\newglossary[slg]{symbols}{sym}{sbl}{List of Symbols}

\usepackage[textsize=small]{todonotes}
\usepackage{varwidth}
\usepackage{url}

\usepackage{enumitem}
\setlist[itemize]{labelindent=*, leftmargin=.5 truecm,nosep}
\usepackage[nottoc]{tocbibind}

\usepackage{fancyhdr}
\usepackage[us,12hr]{datetime} % `us' makes \today behave as usual in TeX/LaTeX
\fancypagestyle{plain}{
\fancyhf{}
\rfoot{ {\ddmmyyyydate\today}}
\lfoot{\thepage}
}

\usepackage{amsthm}

\theoremstyle{plain}
\newtheorem{defi}{Definition}[section]
\newtheorem{prop}[defi]{Proposition}
\newtheorem{theo}[defi]{Theorem}
\newtheorem{theofr}[defi]{Th\'eor\`eme}
\newtheorem{conj}[defi]{Conjecture}
\newtheorem{lemm}[defi]{Lemma}
\newtheorem{coro}[defi]{Corollary}

\theoremstyle{definition}
\newtheorem{rema}[defi]{Remark}
\newtheorem{exem}[defi]{Example}
\newtheorem{exems}[defi]{Examples}

\newcommand{\bdefi}{\begin{defi}}
\newcommand{\edefi}{\end{defi}}
\newcommand{\bprop}{\begin{prop}}
\newcommand{\eprop}{\end{prop}}
\newcommand{\btheo}{\begin{theo}}
\newcommand{\etheo}{\end{theo}}
\newcommand{\btheofr}{\begin{theofr}}
\newcommand{\etheofr}{\end{theofr}}
\newcommand{\blemm}{\begin{lemm}}
\newcommand{\brema}{\begin{rema}}
\newcommand{\erema}{\end{rema}}
\newcommand{\bexer}{\begin{exem}}
\newcommand{\eexer}{\end{exem}}
\newcommand{\bexem}{\begin{exem}}
\newcommand{\eexem}{\end{exem}}
\newcommand{\bexems}{\begin{exems}}
\newcommand{\eexems}{\end{exems}}
\newcommand{\bconj}{\begin{conj}}
\newcommand{\econj}{\end{conj}}
\newcommand{\elemm}{\end{lemm}}
\newcommand{\bcoro}{\begin{coro}}
\newcommand{\ecoro}{\end{coro}}
\newcommand{\dem}{\noindent{\bf Proof. }}

\newcommand{\rem}{\noindent{\bf Remark. }}

\usepackage{mathrsfs}
\renewcommand\mathcal{\mathscr}

\newcommand{\C}{{\cal C}}
\newcommand{\D}{{\cal D}}

\newcommand{\F}{{\cal F}}
\newcommand{\G}{{\cal G}}
\renewcommand{\H}{{\cal H}}

\newcommand{\M}{{\cal M}}

\newcommand{\OOO}{{\cal O}}

\newcommand{\Q}{{\cal Q}}
\newcommand\Rep{\mathcal R}

\newcommand{\Scal}{{\mathcal S}}

\newcommand{\W}{{\cal W}}

\newcommand{\maths}[1]{{\mathbb #1}}  

\newcommand{\CC}{\maths{C}}

\newcommand{\HH}{\maths{H}}

\newcommand{\NN}{\maths{N}}

\newcommand{\PP}{\maths{P}}
\newcommand{\QQ}{\maths{Q}}
\newcommand{\RR}{\maths{R}}

\newcommand{\ZZ}{\maths{Z}}

\newcommand{\mmm}{{\mathfrak m}}

\newcommand{\bs}{\backslash}
\newcommand{\ga}{\gamma}
\newcommand{\Ga}{\Gamma}
\newcommand{\ov}[1]{{\overline{#1}}} 
\newcommand{\ra}{\rightarrow}

\newcommand{\card}{{\operatorname{Card}}}
\newcommand{\CAT}{\operatorname{CAT}}

\newcommand{\covol}{\operatorname{Covol}}

\newcommand{\cqfd}{\hfill$\Box$}

\newcommand{\Det}{\operatorname{Det}}

\newcommand{\gengeod}%{\operatorname{\check{\G}}}
{\operatorname{\widecheck{\G\,}\!\!}}%\overline{\,\!\G}}}

\renewcommand{\Im}{{\operatorname{Im}}}

\newcommand{\n}{\operatorname{\tt n}}

\newcommand{\tr}{\operatorname{\tt tr}}

\newcommand{\Vol}{\operatorname{Vol}}

\newcommand{\hdr}{{\HH}^2_\RR}

\newcommand{\htr}{{\HH}^3_\RR}

\newcommand{\hcr}{{\HH}^5_\RR}

\newcommand{\hnr}{{\HH}^n_\RR}

\newcommand{\hnc}{{\HH}^n_\CC}

\newcommand{\GL}{\operatorname{GL}}

\newcommand{\SL}{\operatorname{SL}}
\newcommand{\PSL}{\operatorname{PSL}}
\newcommand{\SO}{\operatorname{SO}}

\newcommand{\SLO}{\operatorname{SL}_{2}(\OOO)}

\newcommand{\SLH}{\operatorname{SL}_{2}(\HH)}
\newcommand{\PSLH}{\operatorname{PSL}_{2}(\HH)}

\newcommand{\automH}{\operatorname{SU}_f(\OOO)}

\usepackage{interval}
\usepackage{mathtools}
\makeatletter
\DeclareRobustCommand\widecheck[1]{{\mathpalette\@widecheck{#1}}}
\def\@widecheck#1#2{%
    \setbox\z@\hbox{\m@th$#1#2$}%
    \setbox\tw@\hbox{\m@th$#1%
       \widehat{%
          \vrule\@width\z@\@height\ht\z@
          \vrule\@height\z@\@width\wd\z@}$}%
    \dp\tw@-\ht\z@
    \@tempdima\ht\z@ \advance\@tempdima2\ht\tw@ \divide\@tempdima\thr@@
    \setbox\tw@\hbox{%
       \raise\@tempdima\hbox{\scalebox{1}[-1]{\lower\@tempdima\box
\tw@}}}%
    {\ooalign{\box\tw@ \cr \box\z@}}}
\makeatother

\newcounter{fig}

\usepackage{caption}
\usepackage{float}

\title{Integral binary Hamiltonian forms and their waterworlds}
\author{Jouni Parkkonen \and Fr\'ed\'eric Paulin} 
\date{\today}

\begin{document}
\bibliographystyle{../alphanum}
\maketitle
\begin{abstract} 
  We give a graphical theory of integral indefinite binary Hamiltonian
  forms $f$ analogous to the one of Conway for binary quadratic forms
  and the one of Bestvina-Savin for binary Hermitian forms. Given a
  maximal order $\OOO$ in a definite quaternion algebra over $\QQ$, we
  define the {\em waterworld} of $f$, analogous to Conway's {\em
    river} and Bestvina-Savin's {\em ocean}, and use it to give a
  combinatorial description of the values of $f$ on $\OOO\times \OOO$.
  We use an appropriate normalisation of Busemann distances to the
  cusps (with an algebraic description given in an independent
  appendix), and the $\SLO$-equivariant Ford-Voronoi cellulation of
  the real hyperbolic $5$-space.
\end{abstract}

\section{Introduction}
\label{sec:intro}

In the beautiful little book \cite{Conway97} (see also
\cite{Weissman17, Hatcher18}), Conway uses Serre's tree $X_\ZZ$ of the
modular lattice $\SL_2(\ZZ)$ in $\SL_2(\RR)$ (see \cite{Serre83}),
considered as an equivariant deformation retract of the upper
halfplane model of the hyperbolic plane $\hdr$, in order to give a
graphical theory of binary quadratic forms $f$. The components $C$ of
$\hdr-X_\ZZ$ consist of points closer to a given cusp $p/q\in
\PP^1(\QQ)$ of $\SL_2(\ZZ)$ than to all the other ones. When $f$ is
indefinite, anisotropic and integral over $\ZZ$, Conway constructs a
line $R(f)$ in $X_\ZZ$, called the {\it river} of $f$, separating the
components $C$ of $\hdr-X_\ZZ$ such that $f(p,q)>0$ from the ones with
$f(p,q)<0$. This allows a combinatorial description of the values
taken by $f$ on integral points.

Bestvina and Savin in \cite{BesSav12} have given an analogous
construction when $\RR$ is replaced by $\CC$, $\ZZ$ by the ring of
integers $\OOO_K$ of a quadratic imaginary extension $K$ of $\QQ$,
$\hdr$ by $\htr$ and $X_\ZZ$ by Mendoza's spine $X_{\OOO_K}$ in $\htr$
for the Bianchi lattice $\SL_2(\OOO_K)$ in $\SL_2(\CC)$ (see
\cite{Mendoza80}). They construct a subcomplex $O(f)$ of $X_{\OOO_K}$,
called the {\it ocean} of $f$, for any indefinite anisotropic integral
binary Hermitian form $f$ over $\OOO_K$, separating the components of
$\htr-X_{\OOO_K}$ on whose point at infinity $f$ is positive from the
negative ones, and prove that it is homeomorphic to a $2$-plane.

In this paper, we give analogs of these constructions and results for
Hamilton's quaternions and maximal orders in definite quaternion
algebras over $\QQ$.

Let $\HH$ be the standard Hamilton quaternion algebra over $\RR$, with
conjugation $x\mapsto \overline{x}$, reduced norm $\n$ and reduced
trace $\tr$. Let $\OOO$ be a maximal order in a quaternion algebra $A$
over $\QQ$, which is definite (that is, $A\otimes_\QQ \RR=\HH$), with
class number $h_A$ and discriminant $D_A$. An example is given by the
{\it Hurwitz order} $\OOO=\ZZ+\ZZ i+\ZZ j+ \ZZ\frac{1+i+j+k}{2}$, in
which case $h_A=1$ and $D_A=2$. We refer for more information to
\cite{Vigneras80} and Subsection \ref{subsect:quaternion}.  The {\it
  Hamilton-Bianchi group} $\SLO$, which is defined using Dieudonn\'e
determinant, is a lattice in $\SL_2(\HH)$. It acts discretely on the
real hyperbolic $5$-space $\hcr$ with finite volume quotient. The
number of cusps of the hyperbolic orbifold $\SLO\bs\hcr$ is ${h_A}^2$
by \cite[Satz 2.1, 2.2]{KraOse90}, see also \cite[\S 3]{ParPau13ANT}.

Analogously to \cite{Mendoza80} in the complex case, we give in
Section \ref{sec:reduction} an appropriate normalisation of the
Busemann distance to the cusps, and we construct the Ford-Voronoi cell
decomposition of $\hcr$ for $\SLO$, so that the interior of the {\it
  Ford-Voronoi cell} $\H_\alpha$ consists of the points in $\hcr$
closer to a given cusp $\alpha\in \PP^1_r(A)$ of $\SLO$ than to all
the others.  If $X_\OOO$ is the codimension $1$ skeleton of the
Ford-Voronoi cellulation, called the {\it spine} of $\SLO$, then the
hyperbolic $5$-orbifold $\SLO\bs\hcr$ retracts by strong deformations
onto the finite $4$-dimensional orbihedron $\SLO\bs X_\OOO$. 

Using uniform $3$-, $4$- and $5$-polytopes, we give in Example
\ref{ex:gausspine} when $D_A=2$ and in Example \ref{ex:eisenspine}
when $D_A=3$, a complete description of the quotient $\SLO\bs X_\OOO$
and of the link of its vertex. For instance, if $\OOO$ is the Hurwitz
order, then $\SLO\bs X_\OOO$ is obtained by identifying opposite faces
and taking the quotient of any $4$-dimensional cell of $X_\OOO$ by its
stabilizer.  In this case, a $4$-dimensional cell of $X_\OOO$
identifies with the {\it $24$-cell} (the self-dual convex regular
Euclidean $4$-polytope with Schl\"afli symbol $\{3,4,3\}$), and its
stabilizer is isomorphic with an index $2$ subgroup of the Coxeter
group $[3,4,3]$.

Following H.~Weyl \cite{Weyl40}, we will call {\it Hamiltonian form} a
Hermitian form over $\HH$ with anti-involution the conjugation. We
refer to Subsection \ref{subsect:indefbinhamform} and for instance to
\cite{ParPau13ANT} for background. See \cite{ParPau13ANT}
  also for a sharp asymptotic result on the average number of their
  integral representations. Let $f:\HH\times\HH\ra \RR$ be a binary
Hamiltonian form, with
$$
f(u,v)=a\,\n(u)+ \tr(\ov u\, b \,v) +c\,\n(v)\;,
$$ 
which is {\it integral} over $\OOO$ (its {\it coefficients} $a,b,c$
satisfy $a,c\in\ZZ$ and $b\in\OOO$) and indefinite (its {\it
  discriminant} $\Delta(f)=\n(b)-ac$ is positive).  We choose this
definition of integrality for simplicity as in \cite{ParPau13ANT}, in
order to avoid half-integral coefficient in the matrix of the form.
The {\it group of automorphs} of $f$ is the arithmetic lattice
$$
\automH=\{g\in\SLO\;:\;f\circ g=f\}\;.
$$

If $C$ is a Ford-Voronoi cell for $\SLO$, let $F(C)=
\frac{f(a,b)}{\n(\OOO a+\OOO b)}$ where $ab^{-1}\in \PP^1_r(A)$ is the
cusp of $C$. We will say that $C$ is respectively {\it positive,
  negative or flooded} if $F(C)>0$, $F(C)<0$ or $F(C)=0$. Contrarily
to the real and complex cases, there are always flooded Ford-Voronoi
cells, since by taking a $\ZZ$-basis of $\OOO$, the Hamiltonian form
$f$ becomes an integral binary quadratic form over $\ZZ$ with $8\geq
5$ variables, hence always represents $0$. Our countably many flooded
Ford-Voronoi cells are thus the analogues of Conway's two {\it lakes}
for an indefinite isotropic integral binary quadratic form over $\ZZ$.
On the components of $\hdr-X_\ZZ$ along the lakes, Conway proved that
the values of such a form consist in an infinite arithmetic
progression. An analogous result holds in our case, that we only state
when the class number is one in this introduction in order to simplify
the statement (see Proposition \ref{prop:arithprogress} for the
general result.)

\bprop\label{prop:arithprogressintro} If $h_A=1$, given a flooded
Ford-Voronoi cell $C$, there exists a finite set of nonconstant affine
maps $\{\varphi_i:\HH\ra \RR\;:\;i\in F\}$ defined over $\QQ$ such
that the set of values of $F$ on the Ford-Voronoi cells meeting $C$
is $\bigcup_{i\in F}\varphi_i(\OOO)$.  
\eprop

In order to simplify the next statement, assume from now on in this
introduction that the flooded Ford-Voronoi cells are pairwise
disjoint.
%%
%(as for instance $f(u,v)=\n(u)+b\tr(\overline{u}\,v)$ with
%$b\in\NN-\{0,1,2\}$).
%\todo{still the problem about existence here!}
%%
We define the {\it waterworld} $\W(f)$ of $f$ as the subcomplex of the
spine separating positive Ford-Voronoi cells from negative ones, that
is, $\W(f)$ is the union of the cells of $X_\OOO$ contained in (the
boundary of) both a positive and a negative Ford-Voronoi cell. The
{\it coned-off waterworld} $\C\W(f)$ is the union of $\W(f)$ and, for
all cells $\sigma$ of $\W(f)$ contained in a flooded Ford-Voronoi cell
$\H_\alpha$, of the cone with base $\sigma$ and vertex at infinity
$\alpha$. The following result (see Section \ref{sec:waterworld}) in
particular says that $\C\W(f)$ is a piecewise hyperbolic polyhedral
$4$-plane contained in the spine of $\SLO$ except for its ideal cells.

Let $\C(f)$ be the hyperbolic hyperplane of $\hcr$ whose boundary is the projective
set of zeros $\{[u:v]\in\PP^1_r(\HH)\;:\; f(u,v)=0\}$ of $f$.

\btheo\label{theo:mainintro}
The closest point mapping from the coned-off waterworld $\C\W(f)$ to $\C(f)$ 
is an $\automH$-equivariant homeomorphism.
\etheo

Section \ref{sec:background} recalls the necessary information on the
definite quaternion algebras over $\QQ$, the Hamilton-Bianchi groups,
and the binary Hamiltonian forms. Section \ref{sec:reduction} gives
the construction of the normalized Busemann distance to the cusp, and
uses it to give a quantitative reduction theory \`a la Hermite (see
for instance \cite{Borel66b}) for the arithmetic group $\SLO$.  We
describe the Ford-Voronoi cellulation for $\SLO$ and its spine
$X_\OOO$ in Section \ref{sec:spine}. We define the waterworlds and
prove their main properties in Section \ref{sec:waterworld}. The
noncommutativity of $\HH$ and the isotropic property of $f$ require at
various point of this text a different approach than the one in
\cite{BesSav12}.

Recall (see for instance \cite[\S 7]{ParPau13ANT} and Section
\ref{sec:reduction}) that there is a correspondence between positive
definite binary Hamiltonian forms with discriminant $-1$ and the upper
halfspace model of the real hyperbolic $5$-space. In the independent
Appendix \ref{sec:appAalgebraicdistcusp}, we give an algebraic formula
for the Busemann distance of a point $x\in\hcr$ to a cusp $\alpha\in
\PP^1_r(A)$ in terms of the covolume of the $\OOO$-flag associated
with $\alpha$, with respect to the volume of the positive definite
binary Hamiltonian form associated with $x$, analogous to the one of
Mendoza in the complex case.  Furthermore, in the proof of Theorem
\ref{theo:coverhoroball}, we use the upper bound on the minima of
positive definite binary Hamiltonian forms given in \cite{ChePau19}:
If $\ga_2(\OOO)$ is the upper bound, on all such forms $f$ with
discriminant $-1$, of the lower bound of $f(u,v)$ on all nonzero
$(u,v)\in\OOO\times \OOO$, then
\begin{equation}\label{eq:majoconsthermhamil}
\ga_2(\OOO)\leq \sqrt{D_A}\;.
\end{equation}

\medskip\noindent{\small {\it Acknowledgements: } This work was
supported by the French-Finnish CNRS grant PICS \textnumero\,6950.
The second author greatly acknowledges the financial support of
Warwick University for a one month stay, decisive for the writing of
this paper. We warmly thank the referee of a previous version of this
paper for a long list of useful comments which have improved its
content and presentation.}

\section{Backgrounds} 
\label{sec:background}

We refer to \cite{ParPau13ANT} for more informations on the objects
considered in this paper, and we only recall what is strictly needed.

\subsection{Background on definite quaternion algebras over $\QQ$}
\label{subsect:quaternion}

A {\em quaternion algebra} over a field $F$ is a four-dimensional
central simple algebra over $F$. We refer to \cite{Vigneras80} for
generalities on quaternion algebras.  A real quaternion algebra is
isomorphic either to $\M_2(\RR)$ or to Hamilton's quaternion algebra
$\HH$ over $\RR$, with basis elements $1,i,j,k$ as a $\RR$-vector
space, with unit element $1$ and $i^2=j^2=-1$, $ij=-ji=k$. We define
the {\em conjugate} of $x=x_0+x_1i+x_2j+x_3k$ in $\HH$ by
$\overline{x}=x_0-x_1i-x_2j-x_3k$, its {\em reduced trace} by
$\tr(x)=x+\overline{x}$, and its {\em reduced norm} by $\n(x)=
x\,\overline{x}=\overline{x}\,x$. Note that $\n(xy)=\n(x)\n(y)$,
$\tr(\overline{x})=\tr(x)$ and $\tr(xy)=\tr(yx)$ for all $x,y\in\HH$.
For every matrix $X=(x_{i,j})_{1\leq i\leq p,\; 1\leq j\leq q}\in
\M_{p,q}(\HH)$, we denote by $X^*= (\overline{x_{j,i}})_{1\leq i\leq
  q,\;1\leq j\leq p} \in \M_{q,p}(\HH)$ its adjoint matrix. We endow
$\HH$ with the Euclidean norm $x\mapsto \sqrt{\n(x)}$, making the
$\RR$-basis $1,i,j,k$ orthonormal.

Let $A$ be a quaternion algebra over $\QQ$. We say that $A$ is {\em
  definite} (or ramified over $\RR$) if the real quaternion algebra
$A\otimes_\QQ\RR$ is isomorphic to $\HH$, and we then fix an
identification between $A$ and a $\QQ$-subalgebra of $\HH$.  The {\em
  reduced discriminant} $D_A$ of $A$ is the product of the primes
$p\in\NN$ such that the quaternion algebra $A\otimes_\QQ\QQ_p$ over
$\QQ_p$ is a division algebra.  Two definite quaternion algebras over
$\QQ$ are isomorphic if and only if they have the same reduced
discriminant, which can be any product of an odd number of primes (see
\cite[page 74]{Vigneras80}).

A {\em $\ZZ$-lattice} $I$ in $A$ is a finitely generated $\ZZ$-module
generating $A$ as a $\QQ$-vector space.  An {\it order} in $A$ is a
unitary subring $\OOO$ of $A$ which is a $\ZZ$-lattice. In particular,
$A=\QQ\OOO=\OOO\QQ$. Each order of $A$ is contained in a maximal
order. For instance $\OOO=\ZZ+\ZZ i+\ZZ j +\ZZ\frac{1+i+j+k}{2}$ is a
maximal order, called the {\it Hurwitz order}, in $A= \QQ +\QQ i+\QQ j
+\QQ k$ with $D_A=2$. Let $\OOO$ be an order in $A$. The reduced norm
$\n$ and the reduced trace $\tr$ take integral values on $\OOO$. The
invertible elements of $\OOO$ are its elements of reduced norm
$1$. Since $\overline{x}=\tr(x)-x$, any order is invariant under
conjugation.

The {\em left order} $\OOO_\ell(I)$ of a $\ZZ$-lattice $I$ is $\{x\in
A\;:\; xI \subset I\}$. A {\em left fractional ideal} of $\OOO$ is a
$\ZZ$-lattice of $A$ whose left order is $\OOO$. A {\em left ideal} of
$\OOO$ is a left fractional ideal of $\OOO$ contained in $\OOO$.  A
(left) {\it ideal class} of $\OOO$ is an equivalence class of nonzero
left fractional ideals of $\OOO$ for the equivalence relation
$\mmm\sim \mmm'$ if $\mmm'=\mmm c$ for some $c\in A^\times$. The {\it
  class number} $h_A$ of $A$ is the number of ideal classes of a
maximal order $\OOO$ of $A$. It is finite and independent of the
maximal order $\OOO$, and we have $h_A=1$ if and only if
$D_A=2,3,5,7,13$ (see for instance \cite{Vigneras80}).

The {\it reduced norm} $\n(\mmm)$ of a nonzero left  ideal
$\mmm$ of $\OOO$ is the greatest common divisor of the norms of the
nonzero elements of $\mmm$. In particular, $\n(\OOO)=1$.  By
\cite[p.~59]{Reiner75}, we have
\begin{equation}\label{eq:rednormindex}
\n(\mmm)=[\OOO:\mmm]^{\frac{1}{2}}\;.
\end{equation}
The {\it reduced norm} of a nonzero left fractional ideal
$\mmm$ of $\OOO$ is $\frac{\n(c\mmm)}{\n(c)}$ for any $c\in \NN-\{0\}$
such that $c\mmm\subset \OOO$. By Equation \eqref{eq:rednormindex}, if 
$\mmm,\mmm'$ are nonzero left fractional ideals
of $\OOO$ with $\mmm'\subset \mmm$, we have
\begin{equation}\label{eq:rednormindexfrac}
\frac{\n(\mmm')}{\n(\mmm)}=[\mmm:\mmm']^{\frac{1}{2}}\;.
\end{equation}

For $K=\HH$ or $K=A$, we consider $K\times K$ as a right module over
$K$ and we denote by $\PP^1_r(K)= (K\times K-\{0\})/ K^\times$ the
right projective line of $K$, identified as usual with the Alexandrov
compactification $K\cup\{\infty\}$ where $[1:0] =\infty$ and
$[x:y]=xy^{-1}$ if $y\neq 0$.

\subsection{Background on Hamilton-Bianchi groups}
\label{subsect:hambiagroup}

Refering to \cite{Fueter27,Dieudonne43,Aslaksen96}, the {\em
  Dieudonn\'e determinant} $\Det$ is the group morphism from the group
$\operatorname{GL}_2(\HH)$ of invertible $2\times 2$ matrices with
coefficients in $\HH$ to $\RR^*_+$, defined by
$$
\big(\Det\Big(\begin{array}{cc}a& b\\c& d\end{array}\Big)\big)^2\;=
\n(a\,d)+ \n(b\,c) - \tr(a\,\ov c\,d\,\ov b) \;.
$$
If $c\neq 0$, we have (see loc.~cit.) 
\begin{equation}\label{eq:detdieudbis}
\big(\Det\Big(\begin{array}{cc}a& b\\c& d\end{array}\Big)\big)^2\;=
\n(ac^{-1}dc-bc)\;.
\end{equation}
It is invariant under the adjoint map $g\mapsto g^*$. Let $\SLH$ be
the group of $2\times 2$ matrices with coefficients in $\HH$ and
Dieudonn\'e determinant $1$.  We refer for instance to
\cite{Kellerhals03} for more information on $\SLH$.

The group $\SLH$ acts linearly on the left on the right $\HH$-module
$\HH\times\HH$. The projective action of $\SLH$ on $\PP^1_r(\HH)$,
induced by its linear action on $\HH\times \HH$, is the action by
homographies on $\HH\cup\{\infty\}$ defined by
$$
\Big(\begin{array}{cc} a & b \\ c & d\end{array}\Big)\cdot z =
\left\{\begin{array}{ll} (az+b)(cz+d)^{-1} &
{\rm if}\; z\neq \infty,-c^{-1}d \\
ac^{-1} & {\rm if}\; z=\infty, c\neq 0\\
\infty & {\rm otherwise~.}\end{array}\right.
$$

We use the upper halfspace model $\{(z,r)\;:\;z\in\HH, r>0\}$ with
Riemannian metric $ds^2(z,r)= \frac{ds^2_\HH(z)+dr^2}{r^2}$ for the
real hyperbolic space $\hcr$ with dimension $5$. Its space at infinity
$\partial_\infty\hcr$ is hence $\HH\cup\{\infty\}$.  The action of
$\SLH$ by homographies on $\partial_\infty\hcr$ extends to a left
action on $\hcr$ by
\begin{equation}\label{eq:Poincareextension}
 \Big(\begin{array}{cc} a & b \\ c
    & d\end{array}\Big)\cdot (z,r)= 
\Big(\;\frac{(az+b)\,\overline{(cz+d)}+a\,\overline{c}\,r^2}
{\n(cz+d)+r^2\n(c)}, \, \frac r{\n(cz+d)+r^2\n(c)}\,\Big)\;.
\end{equation}
In this way, the group $\PSLH$ is identified with the group of
orientation preserving isometries of $\hcr$. 

For any order $\OOO$ in a definite quaternion algebra $A$
over $\QQ$, the {\it Hamilton-Bianchi group} 
$$
\Ga_\OOO=\SLO= \SLH \cap\M_2(\OOO)
$$ 
is a nonuniform arithmetic lattice in the connected real Lie group
$\SLH$ (see for instance \cite[page 1104]{ParPau10GT} for details). In
particular, the quotient real hyperbolic orbifold $\Ga_\OOO\bs\hcr$
has finite volume. 

\medskip

\rem It would be very interesting to know if the image in
$\PSL_2(\HH)$ of $\SLO$ is commensurable (up to conjugation) to one of
the lattices in $\operatorname{SO}_0(1,5) \simeq\PSL_2(\HH)$ studied
by Vinberg \cite{Vinberg67}, Allcock \cite{Allcock00}, Everitt
\cite{Everitt04}, Ratcliffe-Tschantz \cite{RatTsc04} and others.

\medskip
Recall that the maximal order $\OOO$ is {\em left-Euclidean} if for
all $a,b\in\OOO$ with $b\neq 0$, there exists $c,d\in\OOO$ with
$a=cb+d$ and $\n(d)<\n(b)$, or, equivalently, if for every $\alpha\in
A$, there exists $c\in\OOO$ such that $\n(\alpha-c)<1$. By for
instance \cite[p.~156]{Vigneras80}, $\OOO$ is left-Euclidean if and
only if $D_A\in\{2,3,5\}$.  The following elementary lemma gives a
nice set of generators for $\SLO$.  For us, it will be useful in
Section \ref{sec:spine}.  See also \cite[\S 4]{Speiser32} and \cite[\S
  8]{JohWei99} for the first claim for the Hurwitz order.

\blemm \label{lem:geneSLO}
If $\OOO$ is left-Euclidean, then the group $\SLO$ is generated by
$J=\begin{pmatrix} 0 & 1\\ 1 &0\end{pmatrix}$, $T_w=\begin{pmatrix} 1
& w\\ 0 &1\end{pmatrix}$ for $w\in\OOO$ and $C_{u,v}=\begin{pmatrix} u
& 0\\ 0 &v\end{pmatrix}$ for $u,v\in\OOO^\times$. In particular, the
anti-homography $z\mapsto \overline{z}$ normalizes the action by
homographies of $\SLO$ on $\HH$.  
\elemm

\medskip
\dem The last claim follows from the first one, since $J^{-1}=J$,
$T_w^{-1}=T_{-w}$, $C_{u,v}^{-1}=C_{u^{-1},v^{-1}}$ and for all $z\in
\HH$, we have
$$
\overline{J\cdot \overline{z}}=J\cdot z,\;\;\;
\overline{T_w\cdot \overline{z}}=T_{\overline{w}}\cdot z,\;\;\;
\overline{C_{u,v}\cdot \overline{z}}=
C_{\overline{v}^{-1},\overline{u}^{-1}}\cdot z\;.
$$

Let $G$ be the subgroup of $\SLO$ generated by the matrices $J,T_w,
C_{u,\,v}$ for $w\in\OOO$ and $u,v\in\OOO^\times$ (their Dieudonn\'e
determinant is indeed $1$). Let us prove that any $M=\begin{pmatrix} a
& b\\ c &d\end{pmatrix}\in\SLO$ belongs to $G$, by induction on the
integer $\n(c)$. If $c=0$, then $M=C_{a,\,d}\; T_{a^{-1}b}$ belongs to
$G$. Otherwise, since $\OOO$ is left-Euclidean, there exists
$w,c'\in\OOO$ such that $a=wc+c'$ and $\n(c')<\n(c)$. Hence
$$
M=\begin{pmatrix} 1 & w\\ 0 &1\end{pmatrix}
\begin{pmatrix} 0 & 1\\ 1 &0\end{pmatrix}
\begin{pmatrix} c & d\\ c'
&b-w\,d\end{pmatrix}
$$ 
belongs to $G$ by induction.  
\cqfd

\bcoro\label{cor:geneSLO} 
If $\OOO$ is left-Euclidean, if $\{w_1,w_2,w_3,w_4\}$ is a $\ZZ$-basis
of $\OOO$ and if $S$ is a generating set of the group of units
$\OOO^\times$, then the set
$$
\{J,T_{w_1}, T_{w_2}, T_{w_3}, T_{w_4}\}\cup \{C_{u,v}\;:\;u,v\in S\}
$$ 
is a generating set for  $\SLO$. \cqfd
\ecoro

\medskip
The action by homographies of the group $\Ga_\OOO=\SLO$ preserves the
right projective space $\PP^1_r(A) = A\cup\{\infty\}$, which is the
set of fixed points of the parabolic elements of $\Ga_\OOO$ acting on
$\hcr\cup\partial_\infty \hcr$. In particular, the topological
quotient space $\Ga_\OOO\bs(\hcr\cup\PP^1_r(A))$ is the
compactification of the finite volume hyperbolic orbifold
$\Ga_\OOO\bs\hcr$ by its (finite) space of ends.

\subsection{Background on binary Hamiltonian forms}
\label{subsect:indefbinhamform}

A binary Hamiltonian form $f$ is a map $\HH \times \HH \ra \RR$ with
$$
f(u, v) = a\n(u) + \tr(\ov u\, b\, v) + c\n(v)\,
$$
whose {\it coefficients} $a=a(f)$, $b=b(f)$ and $c=c(f)$ satisfy
$a,c\in\RR$, $b\in\HH$.  Note that $f((u,v)\lambda)=\n(\lambda)f(u,v)$
for all $u,v,\lambda\in\HH$.

The {\it matrix} $M(f)$ of $f$ is the Hermitian matrix
$\Big(\begin{array}{cc}a& b\\ \ov b& c\end{array}\Big)$, so that
$$
f(u,v)= \Big(\begin{array}{c}\!u\!\\ \!v\! \end{array}\Big)^*\;
\Big(\begin{array}{cc}a& b\\\ov b& c\end{array}\Big)\;
\Big(\begin{array}{c} \!u\!\\ \!v\!\end{array}\Big)\,.
$$ 
The {\it discriminant} of $f$ is
$$
\Delta (f) = \n(b)- ac.
$$ 
An easy computation shows that the Dieudonn\'e determinant of $M(f)$
is equal to $|\Delta (f) |$.  A binary Hamiltonian form is {\it
  indefinite} if takes both positive and negative values. It is easy
to check that a form $f$ is indefinite if and only if $\Delta(f)$ is
positive, see \cite[\S 4]{ParPau13ANT}.

The linear action on the left on $\HH\times\HH$ of the group $\SLH$
induces an action on the right on the set of binary Hamiltonian forms
$f$ by precomposition. The matrix of $f\circ g$ is $M(f\circ g)=
g^*\,M(f)\,g$. For every $g\in \SLH$, we have
\begin{equation}\label{eq:invardiscrim}
\Delta(f\circ g)=\Delta(f)\;.
\end{equation}

For every indefinite binary Hamiltonian form $f$, with $a=a(f)$,
$b=b(f)$ and $\Delta=\Delta(f)$, let
$$
\C_\infty(f)=\{[u:v]\in\PP^1_r(\HH)\;:\;f(u,v)=0\}\;.
$$ 
In $\PP^1_r(\HH)=\HH\cup\{\infty\}$, the set $\C_\infty(f)$ is the
$3$-sphere of center $-\frac{b}{a}$ and radius
$\frac{\sqrt{\Delta}}{|a|}$ if $a\neq 0$, and it is the union of
$\{\infty\}$ with the real affine hyperplane $\{z\in\HH\;:\;
\tr(\overline{z}b)+c=0\}$ of $\HH$ otherwise. The values of $f$ are
positive on (the representatives in $\HH\times\HH$ in) one of the two
components of $\PP^1_r(\HH)- \C_\infty(f)$ and negative on the other
one. The set
$$
\C(f)=\{(z,r)\in\HH\times\,]0,+\infty[\;:\;f(z,1)+a\,r^2=0\}
$$
is the ($4$-dimensional) hyperbolic hyperplane in $\hcr$ with boundary
at infinity $\C_\infty(f)$. For every $g\in\SLH$, we have
\begin{equation}\label{eq:antiequiv}
\C_\infty(f\circ g)= g^{-1}\,\C_\infty(f)\;\;\;{\rm and}\;\;\;
\C(f\circ g)= g^{-1}\,\C(f)\;.
\end{equation}

Given an order $\OOO$ in a definite quaternion algebra over $\QQ$, a
binary Hamiltonian form $f$ is {\it integral} over $\OOO$ if its
coefficients belong to $\OOO$.  Note that such a form $f$ takes
integral values on $\OOO\times \OOO$, but the converse might not be
true. The lattice $\Ga_\OOO=\SLO$ of $\SLH$ preserves the set of
indefinite binary Hamiltonian forms $f$ that are integral over
$\OOO$. The stabilizer in $\Ga_\OOO$ of such a form $f$ is its {\it
  group of automorphs}
$$
\automH=\{g\in\Ga_\OOO\;:\;f\circ g=f\}\;.
$$
If $f$ is integral over $\OOO$, then $\automH\bs\C(f)$ is a finite
volume hyperbolic $4$-orbifold, since $\automH$ is arithmetic and by
Borel-Harish-Chandra's theorem (though it might have been known before
this theorem).

\section{On the reduction theory of binary Hamiltonian forms and 
Hamilton-Bianchi lattices}
\label{sec:reduction}

In this section, we study the geometric reduction theory of positive
definite binary Hamiltonian forms, as in Mendoza \cite{Mendoza80} for
the Hermitian case.  The results will be useful in Section
\ref{sec:waterworld}.  We start by recalling the correspondence
between $\hcr$ and positive definite binary Hamiltonian forms with
discriminant $-1$.

Let $\Q$ be the $6$-dimensional real vector space of binary
Hamiltonian forms, and $\Q^+$ its open cone of positive definite ones.
The multiplicative group $\RR^\times_+$ of positive real numbers acts on
$\Q^+$ by multiplication. We will denote by $\RR^\times_+f$ the orbit
of $f$ and by $\PP_+\Q^+$ the quotient space $\Q^+/\RR^\times_+$.  It
identifies with the image of $\Q^+$ in the projective space $\PP(\Q)$
of $\Q$.

Let $\langle\cdot,\cdot\rangle_\Q$ be the symmetric $\RR$-bilinear form
(with signature $(4,2)$) on $\Q$ such that for every $f\in\Q$,
$$
\langle f,f \rangle_\Q=  -  2  \Delta(f)\;.
$$
That is, for all $f,f'\in\Q$, we have
\begin{equation}\label{eq:formprodscalQ}
\langle f,f' \rangle_\Q= a(f)\;c(f')+c(f)\;a(f')-
\tr(\;\overline{b(f)}\;b(f')\,)\;.
\end{equation}
By Equation \eqref{eq:invardiscrim}, we have, for all $f,f'\in\Q$ and
$g\in\SLH$
\begin{equation}\label{eq:invarscalprodQ}
\langle f\circ g,f'\circ g \rangle_\Q=\langle f,f' \rangle_\Q\;.
\end{equation}

Let $\Q^+_1$ the submanifold of $\Q^+$ consisting of the forms with
discriminant $-1$, and let $\Theta:\hcr\ra \Q^+_1$ be the
homeomorphism such that, for every $(z,r)\in\hcr$,
$$
M(\Theta(z,r))=\frac{1}{r}
\begin{pmatrix} \ \,1 & -z\\ -\overline{z} &\n(z)+r^2\end{pmatrix}\;.
$$ 
The fact that this map is well defined and is a homeomorphism
follows by checking that its composition by the canonical projection
$\Q^+\ra\PP_+\Q^+$ is the inverse of the homeomorphism denoted by
$$
\Phi:\RR_+^\times f \mapsto
\Big(-\frac{b(f)}{a(f)}, \frac{\sqrt{-\Delta(f)}}{a(f)}\;\Big)
$$ 
in \cite[Prop.~22]{ParPau13ANT}. By loc.~cit., the map $\Theta$ is
hence (anti-)equivariant under the actions of $\SLH$~: For all
$x\in\hcr$ and $g\in\SLH$, we have
\begin{equation}\label{eq:antiequivTheta}
\Theta(gx)= \Theta(x)\circ g^{-1}\;.
\end{equation}

Let $\OOO$ be a maximal order in a definite quaternion algebra $A$
over $\QQ$.  For every $\alpha\in A$, let
$$
I_\alpha=\OOO \alpha+\OOO\;,
$$
which is a left fractional ideal of $\OOO$. Let $f_\alpha$ be the
binary Hamiltonian form with matrix
$$
M(f_\alpha)=\frac{1}{\n(I_\alpha)}
\begin{pmatrix} \ \;1 & -\alpha\\ 
-\,\overline{\alpha} &\ \n(\alpha)\end{pmatrix}\;.
$$ 
Note that $f_\alpha$ is a positive scalar multiple of the {\it norm
  form} associated with $\alpha$: for all $z\in\HH$,
$$
f_\alpha(u,v)=\big(\overline{u}\;\overline{v}\big)
\;M(f_\alpha)\begin{pmatrix} u \\ v \end{pmatrix}=
\frac{1}{\n(I_\alpha)} \;\n(u-\alpha v)\;.
$$ 
Besides depending on $\alpha$, the form $f_\alpha$ does depend on
the choice of the maximal order $\OOO$. But its homothety class
$\RR^\times f_\alpha$ depends only on $\alpha$.

Let $f_\infty$ be the binary Hamiltonian form whose matrix is
$M(f_\infty) =\begin{pmatrix} 0 & 0 \\ 0 & 1 \end{pmatrix}$, that is,
$f_\infty: (u,v)\mapsto \n(v)$. Note that for every $\alpha\in
\PP^1_r(A)=A\cup\{\infty\}$, the form $f_\alpha$ is nonzero and
degenerate (its discriminant is equal to $0$), and $\RR^\times
f_\alpha$ belongs to the boundary of $\PP_+\Q^+$ in $\PP(\Q)$. The map
$\Phi^{-1}:\hcr \ra \PP(\Q)$ given by $x\mapsto \RR^\times_+\Theta(x)$
extends continuously to a $\SL_2(A)$-(anti-)equivariant homeomorphism
between $\hcr\cup\PP^1_r(A)$ and its image in $\PP(\Q)$ by sending
$\alpha$ to $\RR^\times f_\alpha$ for every $\alpha\in
\PP^1_r(A)$. Proposition \ref{prop:equivfalpha} below makes precise
the scaling factor for the action of $\SL_2(A)$ on the forms
$f_\alpha$ for $\alpha\in \PP^1_r(A)$.  Its proof will use the
following beautiful (and probably well-known) formula.

\blemm \label{lem:beautiful}
For all $g=\begin{pmatrix} a & b\\ c & d \end{pmatrix}\in\SLH$ and
$z,w\in\HH$ such that $g\cdot  z, g\cdot w\neq \infty$, we have
$$
\n(g \cdot  z - g\cdot w)=\frac{1}{\n(cz+d)\n(c w +d)}\;\n(z-w)\;.
$$
\elemm

\dem Since 
$$
\begin{pmatrix} az+b & aw+b\\cz+d &c w +d\end{pmatrix}=
g\begin{pmatrix} z & w\\1 &1\end{pmatrix}
$$ 
and by taking the square of the Dieudonn\'e determinant (see Equation
\eqref{eq:detdieudbis}), we have
\begin{align*}
&\n(g \cdot  z - g\cdot w)=\n((az+b)(cz+d)^{-1}-(aw+b)(c w +d)^{-1})
\\ =\;&\frac{1}{\n(c w +d)}\;\n((az+b)(cz+d)^{-1}(c w +d)-(aw+b))
\\ =\;& \frac{1}{\n(cz+d)\n(c w +d)}\;
\n((az+b)(cz+d)^{-1}(c w +d)(cz+d)-(aw+b)(cz+d))
\\ =\;& \frac{1}{\n(cz+d)\n(c w +d)}\;\n(z-w)\;.\;\;\;\Box
\end{align*}

\bprop\label{prop:equivfalpha} For all $g=\begin{pmatrix} a & b \\ c &
d \end{pmatrix}\in\SL_2(A)$ and $\alpha=[x:y]\in \PP^1_r(A)$, we have
$$
f_{g\cdot\alpha}\circ g= 
\frac{\n(\OOO x +\OOO y)}{\n(\OOO (ax+by) +\OOO (cx+dy))}\;f_\alpha\;.
$$
\eprop

Note that this implies that $f_{g\cdot\alpha}\circ g= f_\alpha$ if
$g\in\SLO$.

\medskip
\dem The result is left to the reader when $\alpha=\infty $ or
$g\cdot\alpha=\infty$, hence we assume that $\alpha,g\cdot\alpha
\neq\infty$. By Lemma \ref{lem:beautiful}, for all $z\in\HH$ such that
$g\cdot z\neq \infty$, we have
\begin{align*}
f_{g\cdot\alpha}\circ g(z,1) & =\n(cz+d)\;f_{g\cdot\alpha}(g\cdot z,1)=
\frac{\n(cz+d)}{\n(I_{g\cdot\alpha})} \;\n(g\cdot z-g\cdot\alpha)\\
& = \frac{1}{\n(I_{g\cdot\alpha})\n(c\alpha+d)}\n(z-\alpha)=
\frac{\n(I_\alpha)}{\n(I_{g\cdot\alpha})\n(c\alpha+d)}\;f_\alpha(z,1)\;.
\end{align*}
The result easily follows. 
\cqfd

\medskip 
For all $\alpha\in\PP^1_r(A)=A\cup\{\infty\}$ and $x\in\hcr$,
let us define the {\it distance from $x$ to the point at infinity
  $\alpha$} by
$$
d_\alpha(x)=\langle f_\alpha,\Theta(x) \rangle_\Q\;.
$$ 
See Appendix \ref{sec:appAalgebraicdistcusp} for an alternate
description of the map $d_\alpha:\hcr\ra\RR$.

The next result gives a few computations and properties of these maps
$d_\alpha$ (which depend on the choice of maximal order $\OOO$).  We
will see afterwards that $\ln d_\alpha$ is an appropriately normalised
Busemann function for the point at infinity $\alpha$.

\bprop\label{prop:propridalpha}
(1) For all $(z,r)\in\hcr$ and $\alpha\in A$, we have
$$%\begin{equation}\label{eq:formuledistcuspun}
d_\alpha(z,r)=\frac{1}{r\n(I_\alpha)} \big(\n(z-\alpha)+r^2\big)\;,
$$%\end{equation}
and $d_\infty(z,r)=\frac{1}{r}$.

\medskip\noindent
(2) For all $x\in\hcr$ and $\alpha=[u:v]\in \PP^1_r(A)$, we have
$$%\begin{equation}\label{eq:formuledistcuspdeux}
d_\alpha(x)=\frac{\Theta(x)(u,v)}{\n(\OOO u + \OOO v)} \;.
$$%\end{equation}

\noindent
(3) For all $g=\begin{pmatrix} a & b \\ c & d \end{pmatrix} 
\in\SL_2(A)$ and $\alpha=[x:y]\in \PP^1_r(A)$, we have
$$
d_{g\cdot\alpha}\circ g= 
\frac{\n(\OOO x +\OOO y)}{\n(\OOO (ax+by) +\OOO (cx+dy))}\;d_\alpha\;.
$$
\eprop

In particular, if $g\in\SLO$ and $\alpha\in \PP^1_r(A)$, then
$d_{g\cdot\alpha}\circ g= d_\alpha$.

\medskip
\dem (1) Since $M(f_\alpha)=\frac{1}{\n(I_\alpha)}
\begin{pmatrix} \;\ 1 & -\alpha\\ 
-\,\overline{\alpha} &\n(\alpha)\end{pmatrix}$ and 
$M(\Theta(z,r))=\frac{1}{r}
\begin{pmatrix} \ \;1 & -z\\ -\overline{z} &\n(z)+r^2\end{pmatrix}$, 
we have, by Equation \eqref{eq:formprodscalQ},
$$
d_\alpha(z,r)=\langle f_\alpha,\Theta(z,r) \rangle_\Q
=\frac{1}{r\,\n(I_\alpha)}
\big((\n(z)+r^2)+\n(\alpha)-\tr(\,\overline{\alpha}\,z)\big)
=\frac{\n(z-\alpha) +r^2}{r\,\n(I_\alpha)}
\;.
$$
The computation of $d_\infty$ is similar and easier.

\medskip\noindent
(2) Let $x=(z,r)\in\hcr$ and $f=\Theta(x)$. If $v\neq 0$, then
$\alpha=uv^{-1}$, and by the definition of $\Theta$ and Assertion (1),
\begin{align*}
\frac{f(u,v)}{\n(\OOO u + \OOO v)} & =
\frac{f(\alpha,1)}{\n(I_\alpha)}=
\frac{1}{\n(I_\alpha)}\big(\,\overline{\alpha}\;1\big)\;M(f)
\begin{pmatrix} \alpha\\1\end{pmatrix}\\ & =
\frac{\n(\alpha)-\overline{\alpha}\,z -\overline{z}\,\alpha
+\n(z)+r^2}{r\,\n(I_\alpha)}
=\frac{\n(z-\alpha) +r^2}{r\,\n(I_\alpha)}=d_\alpha(x)\;.
\end{align*}
Similarly, if $v=0$, then $\frac{f(u,v)}{\n(\OOO u + \OOO v)} =
f(1,0)= \frac{1}{r}=d_\alpha(x)$.

\medskip\noindent
(3) For every $w\in\hcr$, using the (anti-)equivariance
property of $\Theta$, Equation \eqref{eq:invarscalprodQ} and
Proposition \ref{prop:equivfalpha}, we have 
\begin{align*}
d_{g\cdot\alpha}\circ g(w)&=\langle f_{g\cdot\alpha},\Theta(gw) \rangle_\Q
=\langle f_{g\cdot\alpha},\Theta(w)\circ g^{-1} \rangle_\Q
=\langle f_{g\cdot\alpha}\circ g,\Theta(w) \rangle_\Q\\ & =
\frac{\n(\OOO x +\OOO y)}{\n(\OOO (ax+by) +\OOO (cx+dy))}\;
\langle f_{\alpha},\Theta(w) \rangle_\Q\\ & =
\frac{\n(\OOO x +\OOO y)}{\n(\OOO (ax+by) +\OOO (cx+dy))}
\;d_\alpha(w)\;.\;\;\;\Box
\end{align*}

\medskip
Since $\SLO$ is a noncocompact lattice with cofinite volume in $\SLH$
and set of parabolic fixed points at infinity $\PP^1_r(A)$, there
exists (see for instance \cite{Bowditch93}) a $\Ga$-equivariant family
of horoballs in $\hcr$ centered at the points of $\PP^1_r(A)$, with
pairwise disjoint interiors. Since $\SLO\bs\hcr$ may have several
cusps as mentioned in the introduction, there are various choices for
such a family, and we now use the normalized distance to the points of
$\PP^1_r(A)$ in order to define a canonical such family, and we give
consequences on the structure of the orbifold $\SLO\bs\hcr$.

For all $\alpha\in\PP^1_r(A)$ and $s>0$, we define the {\it normalized
  horoball centered at $\alpha$ with radius $s$} as
$$
B_\alpha(s)=\{x\in\hcr\;:\; d_\alpha(x)\leq s\}\;.
$$ 
The terminology is justified by the following result, which proves in
particular that $B_\alpha(s)$ is indeed a (closed) horoball. Recall
that the Busemann function $\beta:\partial_\infty \hcr \times\hcr
\times \hcr\ra\RR$ is defined, with $t\mapsto \xi_t$ any geodesic ray
with point at infinity $\xi\in\partial_\infty\hcr$, by
$$
(\xi,x,y)\mapsto \beta_\xi(x,y)=\lim_{t\ra+\infty} d(x,\xi_t)-d(y,\xi_t)\;.
$$

\bprop\label{prop:horoballs}
 Let $\alpha\in\PP^1_r(A)$ and $s>0$.

\noindent(1) There exists $c_\alpha\in\RR$ such that $\ln d_\alpha(x)=
\beta_\alpha(x,(0,1)) +c_\alpha$ for every $x\in\hcr$.

\smallskip\noindent(2) If $\alpha\in A$, then $B_\alpha(s)$ is the
Euclidean ball of center $\big(\alpha,\frac{s\n(I_\alpha)}{2}\big)$
and radius $\frac{s\n(I_\alpha)}{2}$.  If $\alpha=\infty$, then
$B_\alpha(s)$ is the Euclidean halfspace consisting of all $(z,r)$
with $r\ge \frac 1s$.

\smallskip\noindent(3) For all $g\in\SLO$, we have
$g(B_{\alpha}(s))=B_{g\cdot\alpha}(s)$.  
\eprop

\dem (1) If $\alpha=\infty$, then for every $(z,r)\in\hcr$, we have
$d_\alpha(z,r)=\frac{1}{r}$ and
$$
\beta_\infty((z,r),(0,1))=\beta_\infty((0,r),(0,1))=-\ln r\;,
$$
hence the result holds with $c_\infty =0$.

If $\alpha\in A$, since $\SL_2(A)$ acts transitively on $\PP^1_r(A)$,
let $g=\begin{pmatrix} a & b \\ c & d \end{pmatrix} \in \SL_2(A)$ be
such that $\alpha=g\cdot\infty$. Recall that the Busemann function is
invariant under the diagonal action of $\SL_2(\HH)$ on
$\partial_\infty \hcr \times\hcr \times \hcr$ and is an additive
cocycle in its two variables in $\hcr$. By Proposition
\ref{prop:propridalpha} (3) since $\infty=[1:0]$, we hence have for
every $x\in\hcr$
\begin{align*}
\ln d_\alpha(x) & = \ln d_{g\cdot\infty}(g(g^{-1}x))
= \ln\frac{d_\infty(g^{-1}x)}{\n(\OOO a+\OOO c)}\\ 
& = \beta_\infty(g^{-1}x,(0,1)) -\ln\n(\OOO a+\OOO c) 
=\beta_{g\cdot\infty}(x,g(0,1)) -\ln\n(\OOO a+\OOO c)
\\ & =\beta_{\alpha}(x,(0,1))+ \beta_{\alpha}((0,1),g(0,1)) 
-\ln\n(\OOO a+\OOO c)\;.
\end{align*}
Hence the result holds, and taking $x=(0,1)$, we have by Proposition
\ref{prop:propridalpha} (1)
$$
c_\alpha=\ln \frac{\n(\alpha)+1}{\n(I_\alpha)}\;.
$$

\medskip
\noindent 
(2) If $\alpha\in A$, for every $(z,r)\in\hcr$, by Proposition
\ref{prop:propridalpha} (1), we have $d_\alpha(z,r)\leq s$ if and only
if $\n(z-\alpha)+r^2 \leq s\,r\n(I_\alpha)$, that is, if and only if
$\n(z-\alpha)+ (r-\frac{s\n(I_\alpha)}{2})^2\leq
\big(\frac{s\n(I_\alpha)}{2}\big)^2$. The second claim of Assertion
(2) is immediate. 

\medskip
\noindent 
(3) This follows from Proposition \ref{prop:propridalpha} (3).
\cqfd

\medskip 
The following result extends and generalizes a result for
$D_A=2$ of \cite[\S 5]{Speiser32}.

\btheo\label{theo:coverhoroball} 
Let $\OOO$ be a maximal order in a definite quaternion algebra $A$
over $\QQ$. 

\smallskip\noindent(1) For all distinct
$\alpha,\beta\in\PP^1_r(A)$, the normalized horoballs $B_\alpha(1)$
and $B_\beta(1)$ have disjoint interior. Furthermore, their
intersection is nonempty if and only if $\alpha=\infty$ and
$\beta\in\OOO$, or $\beta=\infty$ and $\alpha\in\OOO$, or
$\alpha,\beta\neq\infty$ and $I_\alpha I_\beta= \OOO(\alpha-\beta)$,
in which case they meet in one and only one point.

\smallskip\noindent(2) We have
$$
\hcr=\bigcup_{\alpha\in\PP^1_r(A)} \;B_\alpha\big(\,\sqrt{D_A}\,\big)\;.
$$
\etheo

Before proving this result, let us make two remarks. 

(i)
Note that $B_0(1)$ and $B_\infty(1)$ intersect (exactly at their
common boundary point $(0,1)$) whatever the definite quaternion
algebra $A$ over $\QQ$ is. Thus the constant $s=1$ in Assertion (1) is
optimal. The family $(B_\alpha(1))_{\alpha\in\PP^1_r(A)}$ is a
(canonical) family of maximal (closed) horoballs centered at the
parabolic fixed points of $\SLO$ with pairwise disjoint
interiors. Since $\SLO$ is a lattice (hence is geometrically finite
with convex hull of its limit set equal to the whole $\hcr$), the
quotient $\SLO\bs(\hcr-\bigcup_{\alpha\in\PP^1_r(A)} B_\alpha(1)\big)$
is compact (see for instance \cite{Bowditch93}).

(ii) 
Assertion (2) is a quantitative version of the standard geometric
reduction theory (see for instance \cite{GarRag70, Borel62,
  Leuzinger95}) for the structure of the arithmetic orbifold
$\SL_2(\OOO)\bs\hcr$. It indeed implies that if $\Rep$ is a finite
subset of $\SL_2(A)$ such that $\Rep\cdot\infty$ is a set of
representatives of $\SL_2(\OOO)\bs \PP^1_r(A)$, and if $\D_\ga$ is a
fundamental domain for the action on $\HH$ of the stabilizer of
$\infty$ in $\ga^{-1}\SL_2(\OOO)\ga$ for every $\ga\in\Rep$, then a
weak fundamental domain for the action of $\SL_2(\OOO)$ on $\hcr$ is
the finite union $\bigcup_{\ga\in\Rep} \ga\Scal_\ga$ where $\Scal_\ga$
is the Siegel set
$$
\Scal_\ga=(\D_\ga\times\;]0,+\infty[\,)\cap
    \ga^{-1}B_{\ga\cdot\infty}(\sqrt{D_A})\;.
$$

\medskip 
\dem 
(1) Note that two horoballs centered at distinct points at infinity,
which are not disjoint but have disjoint interior, meet at one and
only one common boundary point. Hence the last claim of Assertion (1)
follows from the first two ones.

\begin{figure}[H]
\begin{center}
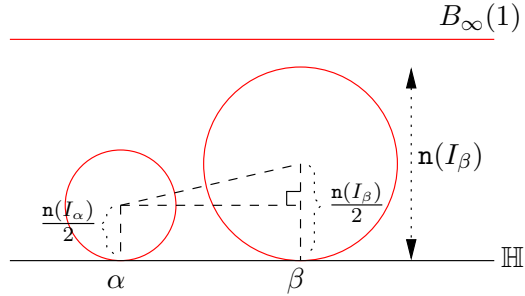
\end{center}
\caption{Disjointness of normalized horoballs $B_{\alpha'}(1)$ for
  $\alpha'\in\PP^1_r(A)$.}\label{fig:disjointhoroballs}
\end{figure}

First assume that $\alpha=\infty$ , so that $\beta\in A$. By
Proposition \ref{prop:horoballs} (2), we have $B_\infty(1)= \{(z,r)
\in\hcr \;:\;r\geq 1\}$ and $B_\beta(1)$ is the horoball centered at
$\beta$ with Euclidean diameter $\n(I_\beta)$ (see Figure
\ref{fig:disjointhoroballs}).  They hence meet if and only if
$\n(I_\beta)\geq 1$, and their interiors meet if and only if
$\n(I_\beta)> 1$. But since $\OOO\subset I_\beta$, by Equation
\eqref{eq:rednormindexfrac}, we have $\n(I_\beta) \leq \n(\OOO)=1$
with equality if and only if $I_\beta=\OOO$, that is, $\beta\in\OOO$.
The result follows.

Up to permuting $\alpha$ and $\beta$ and applying the above argument,
we may now assume that $\alpha,\beta\neq\infty$. The Euclidean balls
$B_\alpha(1)$ and $B_\beta(1)$ meet if and only if the distance
$d_{\alpha\beta}$ between their Euclidean center is less than or equal
to the sum of their radii $r_\alpha$ and $r_\beta$, and their interior
meet if and only if $d_{\alpha\beta}<r_\alpha+r_\beta$. By Proposition
\ref{prop:horoballs} (2) and by the multiplicativity of the reduced
norms (see \cite[Thm.~24.11 and p.~181]{Reiner75}), we have (see the
above picture)
\begin{align*}
{d_{\alpha\beta}}^2 - (r_\alpha+r_\beta)^2 & = \Big(\n(\alpha-\beta)+
\big(\frac{\n(I_\alpha)}{2}-\frac{\n(I_\beta)}{2}\big)^2\Big)-
\big(\frac{\n(I_\alpha)}{2}+\frac{\n(I_\beta)}{2}\big)^2\\ & 
=\n(\alpha-\beta)-\n(I_\alpha)\n(I_\beta)=
\n(\alpha-\beta)-\n(I_\alpha I_\beta)\;.
\end{align*}
Since $\alpha-\beta \in I_\alpha I_\beta$ and again by Equation
\eqref{eq:rednormindexfrac}, we have $\n(\alpha-\beta)\geq \n(I_\alpha
I_\beta)$, with equality if and only if $I_\alpha I_\beta=
\OOO(\alpha-\beta)$. The result follows.

\bigskip
\noindent 
(2) For every $x\in\hcr$, let $(u,v)$ in $\OOO\times\OOO-\{0\}$
realizing the minimum on $\OOO\times\OOO-\{0\}$ of the positive
definite binary Hamiltonian form $\Theta(x)$, whose discriminant is
$-1$. Let $\alpha=[u:v]$. Then by Proposition \ref{prop:propridalpha}
(2) and by Equation \eqref{eq:majoconsthermhamil}, we have, since the
norm of an integral left ideal is at least $1$,
$$
d_\alpha(x)=\frac{\Theta(x)(u,v)}{\n(\OOO u + \OOO v)}\leq \sqrt{D_A}\;.
$$
This proves the result. 
\cqfd 

\medskip The following observation, which is closely related with the
proof of Assertion (1) of Theorem \ref{theo:coverhoroball}, will be
useful later on.

\blemm \label{lem:distBonealphabeta} For all $\alpha\neq\beta\in A$,
the hyperbolic distance between $B_\alpha(1)$ and $B_\beta(1)$ is
$$
d(B_\alpha(1),B_\beta(1))= 
\ln\frac{\n(\alpha-\beta)}{\n(I_\alpha I_\beta)}\;.
$$
\elemm

\dem 
This follows from the easy exercise in real hyperbolic geometry saying
that the distance in the upper halfspace model of the real hyperbolic
$n$-space between two horospheres $\H,\H'$ with Euclidean radius
$r,r'$, and with Euclidean distance between their points at infinity
equal to $\lambda$, is $d(\H,\H')= \ln\frac{\lambda^2}{4rr'}$.

\noindent\begin{minipage}{7.9cm} 
  This exercice uses the facts that the common perpendicular between
  two disjoint horoballs is the geodesic line through their points at
  infinity and that the (signed) hyperbolic length of an arc of
  Euclidean circle centered at a point at infinity with angles with
  the horizontal hyperplane between $\alpha$ and $\pi/2$ is $-\ln \tan
  \frac{\alpha}{2}$. \cqfd\end{minipage} 
\begin{minipage}{7cm}
\begin{figure}[H]
\begin{center}
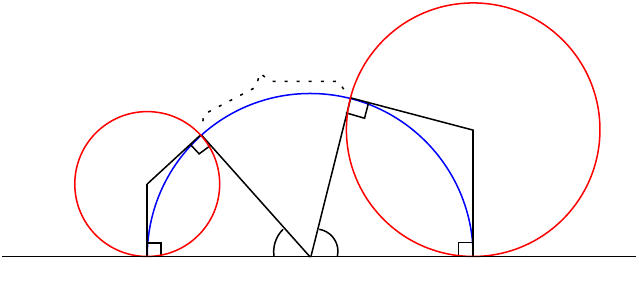
\end{center}
\caption{}
\end{figure}
\end{minipage}

\section{The spine of $\SL_2(\OOO)$}
\label{sec:spine}

Let $A$ be a definite quaternion algebra over $\QQ$ and let $\OOO$ be
a maximal order in $A$. In this section, we describe a canonical
$\SL_2(\OOO)$-invariant cell decomposition of the $5$-dimensional real
hyperbolic space $\hcr$. We follow \cite{Mendoza80,BesSav12} when the
field $\HH$ is replaced by $\CC$, the order $\OOO$ by the ring of
integers of a quadratic imaginary extension of $\QQ$, and $\hcr$ by
$\htr$.

\medskip
For every $\alpha\in\PP^1_r(A)$, the {\em Ford-Voronoi
  cell} of $\alpha$ for the action of $\SL_2(\OOO)$ on $\hcr$
is the set $\H_\alpha$ of points not farther from $\alpha$ than from
any other element of $\PP^1_r(A)$~:
$$
\H_\alpha=\{x\in\hcr\;:\; \forall\;\beta\in\PP^1_r(A),\;\;
d_\alpha(x)\le d_\beta(x)\}\,.
$$ 
In the complex case, this set is called the {\em minimal set} of
$\alpha$, see \cite{Mendoza80}.

\bprop\label{prop:propriFVcells} Let
$\alpha\in\PP^1_r(A)$.

\smallskip\noindent (1) 
For all $g\in\SL_2(\OOO)$, we have $g(\H_\alpha)=\H_{g\cdot\alpha}$.

\smallskip\noindent (2) 
We have $B_\alpha(1)\subset\H_\alpha\subset B_\alpha(\sqrt{D_A})$.

\smallskip\noindent (3)  
The Ford-Voronoi cell $\H_\alpha$ is a noncompact $5$-dimensional
convex hyperbolic polytope, whose proper cells are compact, and the
stabilizer of $\alpha$ in $\SL_2(\OOO)$ acts cocompactly on its
boundary $\partial \H_\alpha$.
 
\smallskip\noindent (4) For every $\beta\in\PP^1_r(A)-\{\alpha\}$, let
  $\Scal_{\alpha,\,\beta}=\{x\in\hcr\;:\; d_\alpha(x)= d_\beta(x)\}$.
  Then $\Scal_{\alpha,\,\beta}$ is a hyperbolic hyperplane, that
  intersects perpendicularly the geodesic line with points at infinity
  $\alpha$ and $\beta$. Furthermore, the Ford-Voronoi cells
  $\H_\alpha$ and $\H_\beta$ have disjoint interior and their
  (possibly empty) intersection is contained in $\Scal_{\alpha,\,\beta}$.

\eprop

Thus
$$
\hcr=\bigcup_{\alpha\in\PP^1_r(A)}\H_\alpha
$$
is a $\SL_2(\OOO)$-invariant cell decomposition of $\hcr$, whose
codimension $1$ skeleton will be studied in the remainder of this
section. We will see in Examples \ref{ex:gausspine} and
\ref{ex:eisenspine} that the inclusions in Assertion (2) of this
proposition, as well as the one of Theorem \ref{theo:coverhoroball},
are sharp when $D_A=2,3$.

\medskip
\dem (1) This follows from Proposition \ref{prop:propridalpha} (3).

\medskip\noindent (2) The inclusion on the left hand side follows from
Theorem \ref{theo:coverhoroball} (1): If $x\in B_\alpha(1)$ and
$x\notin \H_\alpha$, then there exists $\beta\in\PP^1_r(A)-\{\alpha\}$
such that $d_\beta(x)<d_\alpha(x)\leq 1$, thus the interiors of
$B_\alpha(1)$ and $B_\beta(1)$ have nonempty intersection, a
contradiction. If $x\notin B_\alpha(\sqrt{D_A})$, then by Theorem
\ref{theo:coverhoroball} (2), there exists $\beta \in
\PP^1_r(A)-\{\alpha\}$ such that $x\in B_\beta(\sqrt{D_A})$. Hence
$d_\beta(x)\leq \sqrt{D_A} < d_\alpha(x)$, so that $x\notin
\H_\alpha$.

\medskip\noindent (3) and (4) Since $\ln d_\alpha$ is a Busemann
function with respect to the point at infinity $\alpha$ by Proposition
\ref{prop:horoballs} (1), for every $\beta\in\PP^1_r(A)-\{\alpha\}$,
the set $\H_{\alpha,\beta} =\{x\in\hcr\;:\; d_\alpha(x) \le
d_\beta(x)\}$ is a (closed) hyperbolic halfspace. Its boundary is
$\Scal_{\alpha,\,\beta}$, which is hence a hyperbolic hyperplane that
intersects perpendicularly the geodesic line with points at infinity
$\alpha$ and $\beta$. Being the intersection of the locally finite
family of hyperbolic halfspaces $(\H_{\alpha,\beta} )_{\beta\in
  \PP^1_r(A) -\{\alpha\}}$, and containing the horoball $B_\alpha(1)$,
the Ford-Voronoi cell $\H_\alpha$ is a noncompact $5$-dimensional
convex hyperbolic polytope. Since $\alpha$ is a bounded parabolic
fixed point of the lattice $\SL_2(\OOO)$ and by Assertion (2), the
stabilizer of $\alpha$ in $\SL_2(\OOO)$ acts cocompactly on $\partial
\H_\alpha$, and hence the boundary cells of $\H_\alpha$ are compact.
\cqfd

%We can construct a fundamental domain for the action of $\SLO$ in
%$\hcr$ as follows: Let $\ga_i\in\SL_2(A)$ $1\le i\le h_A$ such that
%$\{\gamma_1\infty:1\le i\le h_A\}$ is a set of representatives for the
%$\SLO$-orbits in $A\cup\{\infty\}$, and let $\ga_i\Ga_i\ga_i^{-1}$ be
%the stabiliser of $\alpha_i$ in $\SLO$.  The group $\Ga_i$ has a
%compact fundamental domain $D_i$ for its action in $\HH$. The set
%$$
%F=\bigcup_{i=1}^{h_A}\H_{\ga_i\infty}\cap\ga_i(D_i\times
%\interval[open left]0\infty)
%$$ 
%is a fundamental domain of $\SLO$. The truncation
%$F-\bigcup_{i=1}^{h_A}B_{\ga_i\infty}(1)$ is compact.  See
%\cite[p.~88]{ParPau13ANT} for more details.\footnote{See
% \cite{Siegel65} and \cite[p.~14-15]{Mendoza80}.}

\medskip 
The horoballs $B_0(1)$ and $B_\infty(1)$ with disjoint interiors meet
at $(0,1)\in\hcr$, and at most two horoballs with disjoint interior
can meet at a given point of $\hcr$. Thus, the Ford-Voronoi cells at
$0$ and at $\infty$ have nonempty intersection, which is a compact
$4$-dimensional hyperbolic polytope.  This intersection
$$\Sigma_\OOO=\H_0\cap\H_\infty$$ is called the {\em fundamental cell}
of the spine of $\SL_2(\OOO)$. We will describe it in Example
\ref{ex:gausspine} when $D_A=2$ and in Example \ref{ex:eisenspine}
when $D_A=3$ .

\blemm \label{lem:3dimcell}
Let $\alpha\in\PP^1_r(A)$ be such that $e=\H_\infty\cap \H_0 \cap
\H_\alpha$ is a $3$-dimensional cell in the boundary of $\Sigma_\OOO$.
Then
$$
\min\{\n(I_\alpha),\n(I_{\alpha^{-1}})\} \ge\frac{1}{D_A}\;,
$$ 
and the horizontal projection of $e$ to $\HH$ is contained in the
Euclidean hyperplane 
$$
\{z\in\HH\;:\; \tr(\ov\alpha\, z)=1+\n(\alpha)-\n(I_\alpha)\}\,.
$$
\elemm

\dem %The proof is similar to the one of \cite[Lemma 4.2]{BesSav12}.
Note that $\alpha\neq 0,\infty$. By Proposition
\ref{prop:propriFVcells} (2), the intersection
$B_\infty(\sqrt{D_A})\cap B_0(\sqrt{D_A})\cap B_\alpha(\sqrt{D_A})$
contains $e$, hence the intersections $B_\infty(\sqrt{D_A})\cap
B_\alpha(\sqrt{D_A})$ and $B_0(\sqrt{D_A})\cap B_\alpha(\sqrt{D_A})$
are nonempty. Since $B_\infty(\sqrt{D_A})$ is the Euclidean halfspace
of points $(z,r)$ with $r\geq \frac{1}{\sqrt{D_A}}$ and
$B_\alpha(\sqrt{D_A})$ is a Euclidean ball tangent to the horizontal
plane with diameter $\sqrt{D_A}\n(I_\alpha)$ by Proposition
\ref{prop:horoballs} (2), this implies that $\sqrt{D_A}\n(I_\alpha)
\geq \frac{1}{\sqrt{D_A}}$, so that $D_A\n(I_\alpha)\geq 1$.  Since
$g=\begin{pmatrix} 0 & 1\\ 1 & 0\end{pmatrix}$ belongs to
$\SL_2(\OOO)$ and maps $0$ to $\infty$ and $\alpha$ to $\alpha^{-1}$,
and by Proposition \ref{prop:horoballs} (3), the intersection
$B_\infty(\sqrt{D_A})\cap B_{\alpha^{-1}}(\sqrt{D_A})$ is nonempty,
hence similarly $D_A\n(I_{\alpha^{-1}})\geq 1$.

\begin{figure}[H]
\begin{center}
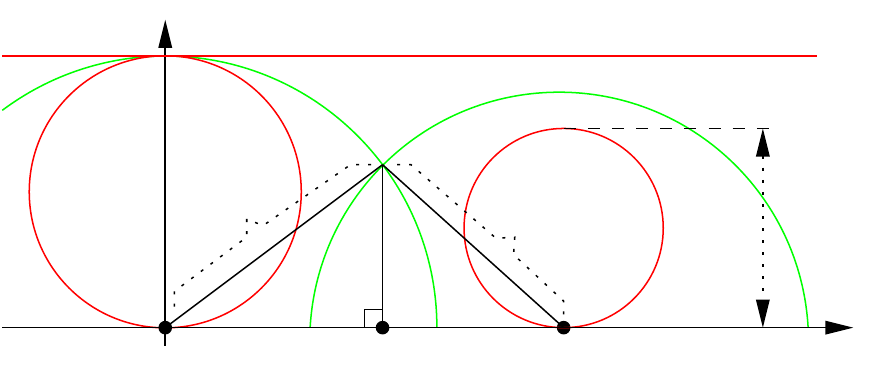
\end{center}
\caption{}\label{fig:case1}
\end{figure}

The set of points equidistant to $0$ and $\infty$ is the open
Euclidean upper hemisphere of radius $1$ centered at $0$, and the set
of points equidistant to $\alpha$ and $\infty$ is the open Euclidean
upper hemisphere of radius $\sqrt{\n(I_\alpha)}$ centered at
$\alpha$. The projection to $\HH$ of the intersection of these
hemispheres is contained in the affine Euclidean hyperplane of $\HH$
perpendicular to the real vector line containing $\alpha$ that passes
through the projection $\lambda \alpha$ with $\lambda>0$ to that line
of any point at Euclidean distance $1$ from $0$ and at Euclidean
distance $\sqrt{\n(I_\alpha)}$ from $\alpha$. An easy computation
(considering the two cases when $\n(\alpha)>1$ as in Figure \ref{fig:case1}
or when $\n(\alpha)\leq 1$) using right angled triangles gives that
$\lambda =\frac{1+\n(\alpha)- \n(I_\alpha)}{2\n(\alpha)}$. Since
$(u,v)\mapsto \frac{1}2 \tr(\ov{u}\,v)$ is the standard Euclidean
scalar product on $\HH$, this gives the result.
\cqfd

\medskip
The {\em spine} of $\SL_2(\OOO)$ is the codimension $1$ skeleton of
the cell decomposition into Ford-Voronoi cells of $\hcr$, that is
$$
X_\OOO=\bigcup_{\alpha\ne\beta\in\PP^1_r(A)}\H_\alpha\cap\H_\beta\;\;=
\bigcup_{\alpha\in\PP^1_r(A)}\partial \H_\alpha\;.
$$ 
It is an $\SL_2(\OOO)$-invariant piecewise hyperbolic polyhedral
complex of dimension $4$. We refer for instance to \cite{BriHae99} for
the definitions related to polyhedral complexes, $\CAT(0)$ spaces and
orbihedra.  Note that the stabilizers in $\SL_2(\OOO)$ of the cells of
$X_\OOO$ may be nontrivial. The spine is called the {\em minimal
  incidence set} in the complex case in \cite{Mendoza80} and \cite{SchVog83}, and the
{\em cut locus of the cusp} in \cite[\S 5]{HerPau02a} when the class
number is one.

For every hyperbolic cell $C$ of $X_\OOO$ and every $\alpha\in
\PP^1_r(A)$ such that $C\subset \partial \H_\alpha$, the radial
projection along geodesic rays with point at infinity $\alpha$ from
$C$ to the horosphere $\partial B_\alpha(1)$ is a homeomorphism onto
its image, and the pull-back of the flat induced length metric on this
horosphere endows $C$ with a structure of a compact Euclidean
polytope. This Euclidean structure does not depend on the choice of
$\alpha$, since the (possibly empty) intersection $\H_\alpha\cap
\H_\beta$ is equidistant to $B_\alpha(1)$ and $B_\beta(1)$ for all
distinct $\alpha,\beta$ in $\PP^1_r(A)$. It is well known (see for
instance \cite{Aitchison10}) that these Euclidean structures on the
cells of $X_\OOO$ endow $X_\OOO$ with the structure of a $\CAT(0)$
piecewise Euclidean polyhedral complex.

Furthermore, $X_\OOO$ is a $\SLO$-invariant deformation retract of
$\hcr$ along the geodesic rays with points at infinity the points in
$\PP^1_r(A)$ and since the quotient orbifold with boundary
$\SL_2(\OOO) \bs \big( \hcr- \bigcup_{\alpha\in\PP^1_r(A)}
B_\alpha(1)\big)$ is compact, the quotient space $\SLO\bs X_\OOO$ is a
finite locally $\CAT(0)$ piecewise Euclidean orbihedral complex.

\medskip 
%In the remainder of  this section, we  describe the cell structure of $\SLO\bs
%X_\OOO$ in some particular cases. 
The following result gives a description of the cell structure of $\SLO\bs
X_\OOO$ when $\OOO$ is left-Euclidean. 
See Examples \ref{ex:gausspine} and \ref{ex:eisenspine} for a more
detailed study when $D_A=2,3$.

\bprop \label{prop:euclideancase}
The Hamilton-Bianchi group $\SLO$ acts transitively on the set
of $4$-dimen\-sional cells of its spine $X_\OOO$ if and only if $D_A\in
\{2,3,5\}$. In these cases, the horizontal projection of the fundamental
cell $\Sigma_\OOO$ to $\HH$ is the Euclidean Voronoi cell  of $0$ for
the $\ZZ$-lattice $\OOO$ in the Euclidean space $\HH$.  
\eprop

\dem 
If $\SLO$ acts transitively on the $4$-dimensional cells of $X_\OOO$,
then $X_\OOO = \SLO\,\Sigma_\OOO$, and the stabilizer of $\infty$ in
$\SLO$ acts transitively on the set of $4$-dimensional cells in
$\partial \H_\infty$, since $\begin{pmatrix} 0 & 1 \\ 1 &
  0 \end{pmatrix} \in\SL_2(\OOO)$ preserves $\Sigma_\OOO = \H_\infty
\cap \H_0$ and exchanges $\H_\infty$ and $\H_0$.  This stabilizer
consists of the upper triangular matrices with coefficients in $\OOO$,
hence with diagonal coefficients in $\OOO^\times$. The orbit of
$0\in\HH$ under this stabilizer is exactly $\OOO$.  Since
$\Sigma_\OOO$ is compact and contained in the open Euclidean upper
hemisphere centered at $0$ with radius $1$, by horizontal projection
on $\HH$, this proves that $\HH$ is covered by the open balls of
radius $1$ centered at the points of $\OOO$. Hence $\OOO$ is
left-Euclidean.

Conversely, if $\OOO$ is left-Euclidean, then the class number of $A$
is $1$, and $\SLO$ acts transitively on the Ford-Voronoi cells. In
order to prove that $\SLO$ acts transitively on the $4$-dimensional
cells of $X_\OOO$, we hence only have to prove that the stabilizer of
$\infty$ in $\SLO$ acts transitively on the $4$-dimensional cells of
$\partial \H_\infty$. For this, let $\alpha\in A$ be such that
$\H_\infty \cap\H_\alpha$ is a $4$-dimensional cell in $\partial
\H_\infty$. Let us prove that $\alpha\in\OOO$, which gives the result.
Due to problems caused by the noncommutativity of $\HH$, the proof of
\cite[Prop.~4.3]{BesSav12} does not seem to extend exactly. We will
use instead Lemma \ref{lem:geneSLO}.

\begin{figure}[H]
\begin{center}
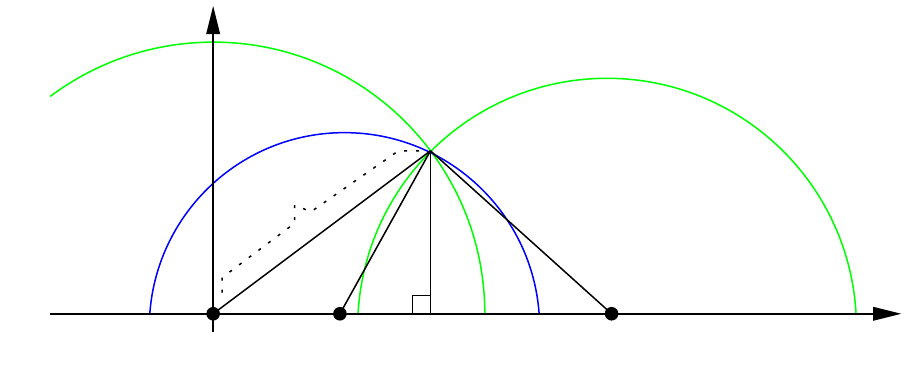
\end{center}
\caption{}\label{fig:inversion}
\end{figure}

Assume for a contradiction that $\alpha\notin \OOO$. Since $\OOO$ is
left-Euclidean, there exists $c\in\OOO$ such that $\n(\alpha-c)<1$.
Up to replacing $\alpha$ by $\alpha-c$, since translations by $\OOO$
preserve $\H_\infty$, we may assume that $0< \n(\alpha) <1$. For every
$\beta\in A$ and $\beta'\in\PP^1_r(A)-\{\beta\}$, let us denote by
$S_{\beta,\beta'}$ the Euclidean upper hemisphere centered at $\beta$
equidistant from the points at infinity $\beta$ and $\beta'$. In
particular, $S_{0,\infty}$ has radius $1$. The inversion with respect
to the sphere containing $S_{0,\infty}$ acts by an
orientation-reversing isometry on $\hcr$, and acts on the boundary at
infinity $\PP^1_r(\HH) =\HH\cup\{\infty\}$ by $z\mapsto \frac{z}
{\n(z)} = \frac{1} {\overline{z}}$. By Lemma \ref{lem:geneSLO}, it
hence normalizes $\SLO$ and, in particular, sends $S_{\alpha,\infty}$
to $S_{\frac{\alpha}{\n(\alpha)},0}$, and fixes $S_{0,\infty}$ (see
Figure \ref{fig:inversion}).  Since $\n(\alpha)<1$, the hemisphere
$S_{\alpha,\infty}$ is therefore below the union of $S_{0,\infty}$ and
$S_{\frac{\alpha}{\n(\alpha)},0}$, which contradicts the fact that
$\H_\infty\cap\H_\alpha$, which is contained in $S_{\alpha,\infty}$,
is a $4$-dimensional cell in $\partial \H_\infty$.

\medskip 
In order to prove the last claim of Proposition
\ref{prop:euclideancase}, note that $\n(I_\alpha)=1$ if $\alpha
\in\OOO$, and that the above proof shows that the $4$-dimensional
cells contained in $\partial\H_\infty$ and meeting the fundamental
cell along a $3$-dimensional cell are contained in spheres centered at
points in $\OOO$. Therefore, by Lemma \ref{lem:3dimcell}, the
horizontal projection of $\Sigma_\OOO$ is the intersection of the
halfspaces containing $0$ and bounded by the Euclidean hyperplanes
with equation $\tr(\overline{\alpha}\, z)=\n(\alpha)$ for all
$\alpha\in\OOO$. Since this hyperplane is the set of points $z$ in the
Euclidean space $\HH$ equidistant to $0$ and $\alpha$, this proves
that the horizontal projection of $\Sigma_\OOO$ is indeed the Voronoi
cell at $0$ of the $\ZZ$-lattice $\OOO$.
\cqfd

\bexem\label{ex:gausspine} Let $A=\QQ+\QQ i +\QQ j +\QQ k\subset \HH$
be the definite quaternion algebra over $\QQ$ with $D_A=2$, and let
$\OOO=\ZZ+\ZZ i +\ZZ j +\ZZ \frac{1+i+j+k}{2}$ be the (maximal)
Hurwitz order in $A$.  The Hurwitz order $\OOO$ is the lattice of type
$F_4=D_4^*$. The group of unit Hurwitz quaternions has
elements  $$\OOO^\times=\Big\{\pm 1,\pm i, \pm j,\pm k, \frac{\pm 1
  \pm i \pm j \pm k}{2}\Big\}\,.$$ 

The Voronoi cell $\Sigma_\OOO^\HH$ of $0$ for the lattice $\OOO$ in
$\HH$ is (up to homothety) the {\it $24$-cell}, which is the (unique)
self-dual, regular, convex Euclidean $4$-polytope, whose Schl\"afli
symbol is $\{3,4,3\}$.  The vertices of $\Sigma_\OOO^\HH$ are the $24$
quaternions
$$
\frac{1+i}2\OOO^\times=
\Big\{ \frac{\pm 1\pm i}2,\frac{\pm 1\pm j}2,\frac{\pm 1\pm k}2,
\frac{\pm i\pm j}2,\frac{\pm i\pm k}2,\frac{\pm j\pm k}2
\Big\}\,.
$$ 
See for instance \cite[p.~119]{ConSlo88} for more details and
references.

Let $\HH^\times_1$ be the subgroup of $\HH^\times$ that consists of
the quaternions of norm $1$. The group morphism $\HH^\times_1 \times
\HH^\times_1\to\SO(4)$ that associates to $(u,v)\in \HH^\times_1\times
\HH^\times_1$ the orthogonal transformation $z\mapsto uzv^{-1}$ is
surjective, see for instance \cite[Thm.~8.9.8]{Berger79}.  The group
of Euclidean symmetries of the $24$-cell is the Coxeter group
$[3,4,3]$. It consists of the $1152$ elements $z\mapsto uzv^{-1}$,
$z\mapsto u\,\ov z\,v^{-1}$ of $\operatorname{O}(4)$, where either $u$
and $v$ are unit Hurwitz integers or $u/\sqrt 2$ and $v/\sqrt 2$ are
in $\frac{1+i}2\OOO^\times$.

By Proposition \ref{prop:euclideancase}, the
fundamental cell of $\SLO$ is
$$
\Sigma_\OOO=\{(z,t)\in\hcr:z\in \Sigma_\OOO^\HH,\ \n(z)+t^2=1\}.
$$
With the notation of Lemma \ref{lem:geneSLO}, the stabilizer of
$\Sigma_\OOO$ in $\SLO$ consists of the $1152$ matrices $C_{a,d}=
\begin{pmatrix} a&0\\0&d\end{pmatrix}$ and $JC_{a,d}=\begin{pmatrix}
0&d\\a&0\end{pmatrix}$ with $a,d\in\OOO^\times$. When $\Sigma_\OOO$ is 
identified with $\Sigma_\OOO^\HH$ by the horizontal projection, the
diagonal matrices induce by Equation \eqref{eq:Poincareextension}
$288$ rotational symmetries of $\Sigma_\OOO^\HH$, and the antidiagonal
ones induce another $288$ orientation-reversing symmetries, together
forming a subgroup of index $2$ in the Coxeter group $[3,4,3]$.

The quotient $\SLO\bs X_\OOO$ is obtained by identifying the opposite
$3$-dimensional cells of $\Sigma_\OOO$ (which are $24$ regular
octahedra) by translations by elements of $\OOO$, and by forming the
quotient by the stabilizer of $\Sigma_\OOO$.  In particular, all the
vertices of $X_\OOO$ are in the same orbit under $\SLO$.

\begin{figure}[H]
\begin{center}
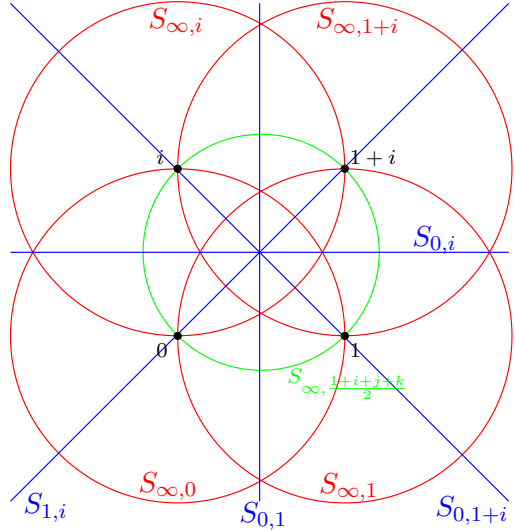
\end{center}
\caption{Boundary of equidistant hemispheres and halfplanes in
$\CC\subset \HH$.}\label{fig:equid}
\end{figure}

Speiser \cite[\S 5]{Speiser32} observed that the estimate of
Proposition \ref{prop:propriFVcells} (2) is sharp in this example:
$\hcr$ is indeed completely contained in $\bigcup_{\alpha\in
  \PP^1_r(A)} B_\alpha(\sqrt 2)$, and the orbit that contains all the
vertices of $\Sigma_\OOO$ is not contained in the union of the
interiors of the horoballs $B_\alpha(\sqrt 2)$.  Furthermore, Speiser
proved that the point
$$
v_0= \Big(\,\frac{1+i}2\,,\frac1{\sqrt 2}\,\Big)
$$
belongs to the boundary of exactly $10$ horoballs $B_\alpha(\sqrt 2)$,
the ones with $\alpha$ in
$$
E=\Big\{\infty, 0,1,i,1+i, \frac{1+i\pm j\pm  k}2, 
\frac 1{1-i}=\frac{1+i}2\Big\}\;.
$$
In particular, $v_0$ is a vertex of the spine $X_\OOO$, contained in
the boundary of exactly $10$ Ford-Voronoi cells $\H_\alpha$ for
$\alpha$ in this set.

The set $E$ contains exactly $5$ pairs $\{\alpha,\beta\}$ of distinct
elements such that the interiors of the horoballs $B_\alpha(\sqrt 2)$
and $B_\beta(\sqrt 2)$ are disjoint, these pairs being $\{\infty,\frac
1{1-i}\}$, $\{0,1+i\}$, $\{1,i\}$, $\{\frac{1+i+j+k}2,
\frac{1+i-j-k}2\}$ and $\{\frac{1+i+j-k}2,\frac{1+i-j+k}2\}$. If
$\{\alpha,\beta\}$ is one of these pairs, the Ford-Voronoi cells
$\H_\alpha$ and $\H_\beta$ intersect only at $v_0$. For all other
pairs in $E$, the intersection is a higher-dimensional cell.

\begin{figure}[H]
\begin{center}
\includegraphics[width=7.3cm]{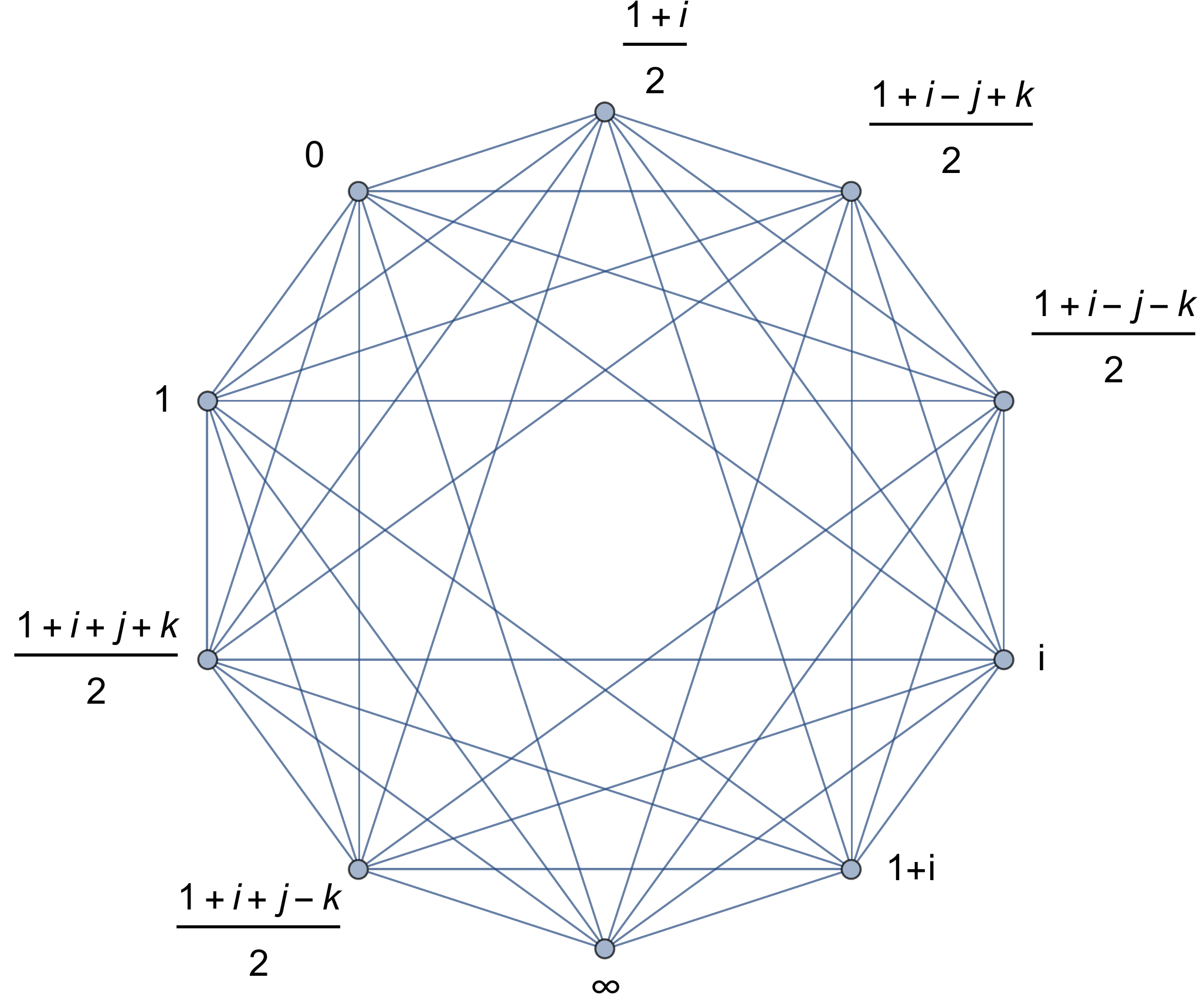}
\end{center}
\caption{Pattern of intersections into $4$-dimensional cells  of
  Ford-Voronoi cells centered at $\{\infty,0,1,i,i+1,\frac{1+i\pm j\pm
    k}{2},\frac{1+i}{2}\}$.}\label{fig:intpatt}
\end{figure}

As $0,1,i,1+i, \frac{1+i\pm j\pm k}2$ are in $\OOO$ and $\frac 1{1-i}$
is not in $\OOO$, there are $8$ Ford-Voronoi cells incident to $v_0$
that intersect $\H_\infty$ in a $4$-dimensional $24$-cell (see Figure
\ref{fig:equid}, which represents the intersection with the plane in
$\HH$ containing $0,1,i$ of the closures of the equidistant spheres
and planes between some pairs of elements in $\{\infty, 0,1,i,
1+i,\frac{1+i+ j+ k}2\}$, so that the horizontal projection of $v_0$
is the common intersection points of the straight lines).  A similar
property holds for all the other Ford-Voronoi cells incident to $v_0$:
For example, $\H_0$ intersects in a $4$-dimensional cell the
Ford-Voronoi cells $\H_\infty,\H_1,\H_i, \H_{\frac{1+i\pm j\pm k}{2}},
\H_{\frac{1+i}{2}}$, but not $\H_{i+1}$ by Theorem
\ref{theo:coverhoroball} (1), since $I_0I_{i+1}=\OOO\ne \OOO(1+i)$.
Thus the pattern of pairwise intersections into $4$-dimensional cells
of these $10$ Ford-Voronoi cells is given by Figure \ref{fig:intpatt}
and the number of $24$-cells containing $v_0$ is exactly $40=(10\times
8)/2$, one for each edge of this intersection pattern.

The boundary of each $\H_\alpha$ is tiled by $24$-cells,
combinatorially forming the $24$-cell honeycomb.  The dual of this
honeycomb is the $16$-cell honeycomb. Therefore, the link
%\footnote{The complex obtained as the intersection of $\partial
%\H_\alpha$ and a sphere with small radius centered at $v_0$.} 
of the vertex $v_0$ in the tessellation of $\partial\H_\alpha$ for all
$\alpha\in E$ is the dual of the $16$-cell, which is the boundary of
the $4$-cube, such that the intersection of the link with each of the
eight $24$-cells is a $3$-cube.
%A small neighbourhood in $\partial\H_\alpha$ of the vertex $v_0$ is
%the cone over the $3$-skeleton of the $4$-cube.

Gluing together the ten boundaries of $4$-cubes (that have been
subdivided in eight $3$-cubes each) according to the above
intersection pattern proves that the link of $v_0$ in the spine
$X_\OOO$ is the $3$-skeleton of the $5$-cube (which is the
$5$-dimensional regular polytope with Schl\"afli symbol $\{4,3,3,3\}$).
\eexem

\bexem \label{ex:eisenspine}
The maximal order of the definite quaternion algebra
$\big(\frac{-1,-3}{\QQ}\big)$ of discriminant $D_A=3$ is $\ZZ[1, i,
  \frac{i + j}2, \frac{1 + k}2]$, see \cite[p.~98]{Vigneras80}.  Using 
the unique $\QQ$-linear map from $\big(\frac{-1,-3}{\QQ}\big)$ to
$\HH$ sending $1$ to $1$, $i$ to $j$, $j$ to $k\sqrt{3}$ and $k$ to
$-i\sqrt{3}$, we identify $\big(\frac{-1,-3}{\QQ}\big)$ with the
$\QQ$-subalgebra $A$ of $\HH$ generated by $1$, $i\sqrt 3 $, $j$ and
$k\sqrt 3 $, and the maximal order is then identified with
$\OOO=\ZZ[1, \rho, j,\rho j]$, where
$$
\rho=\frac{1 + i\sqrt 3 }2\;.
$$  
The group of units of $\OOO$ is the dicyclic group of  order $12$
$$
\OOO^\times=
\{\pm 1,\;\pm j,\;\pm \rho,\;\pm\rho^2,\;\pm \rho j,\; \pm\rho^2 j\}\;.
$$

\begin{figure}[H]
\begin{center}
\includegraphics[width=6cm]{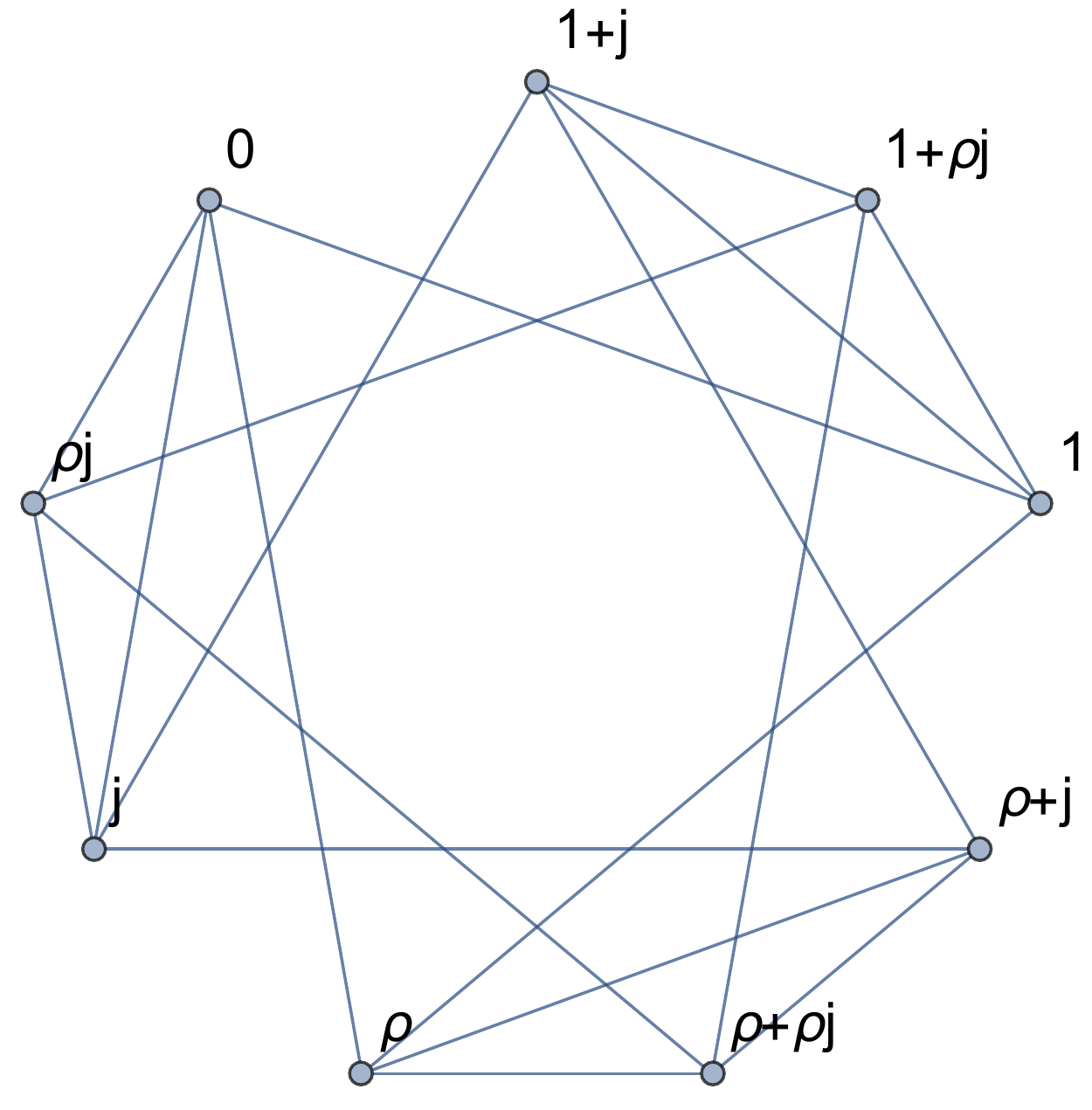}
\end{center} 
\caption{The $1$-skeleton of the $3\,$-$3$
  duoprism.}\label{fig:33duoprismskeleton}
\end{figure}

The elements of the maximal order $\OOO=\ZZ[1, \rho]+\ZZ[1, \rho]j$ of
$A$ are the vertices of the $3\,$-$3$ duoprism honeycomb in the
$4$-dimensional Euclidean space $\HH$. The $9$ elements of the set
$$
V_{3,3}=\{0,1,j,1+j,\rho, \rho j, 1+\rho j,j+\rho,
\rho(1+j)\}\;,
$$ 
contained in $\OOO$, are the vertices of its fundamental $3\,$-$3$
{\em duoprism}, which is a uniform $4$-polytope with Schl\"afli symbol
$\{3\}\times\{3\}$ (the Cartesian product of two equilateral
triangles, whose $1$-skeleton is given in Figure
\ref{fig:33duoprismskeleton}).  We refer to Coxeter's three papers
\cite{Coxeter40,Coxeter85, Coxeter88} for notation and references
about uniform polytopes and their Coxeter groups, with the help of the
numerous and beautiful Wikipedia articles.

The Voronoi cell $\Sigma_\OOO^\HH$ of $0$ for the lattice $\OOO$ in
$\HH$ is the $6\,$-$6$ {\em duoprism} whose Schl\"afli symbol is
$\{6\}\times\{6\}$.  It is the Cartesian product of two copies of the
Voronoi cell of $0$ for the hexagonal lattice of the Eisenstein
integers in $\CC$ whose set of vertices is $V_6= \{\pm
\frac{i}{\sqrt{3}}, \pm\frac{1}{2} \pm\frac{i}{2\sqrt{3}}\}$. Thus,
the set of vertices of $\Sigma_\OOO^\HH$ is $V_6+jV_6$. These $36$
vertices, including
$$
z_0=\frac 12+\frac i{2\sqrt 3}+\frac j2+\frac k{2\sqrt 3}=
(j+\rho)(1+\rho)^{-1}\;,
$$
belong to $A$ and all have reduced norm $\frac 23$.

By Proposition \ref{prop:euclideancase}, the fundamental cell
$\Sigma_\OOO$ of $\SLO$ is the subset of the Euclidean unit sphere in
$\hnc$ whose horizontal projection is $\Sigma_\OOO^\HH$.  In
particular, all the vertices of $\Sigma_\OOO$ have Euclidean height
$\frac{1}{\sqrt{3}}$.  Let $u$ and $v$ be either both in $\OOO^\times$
or both in $\rho^{\frac{1}{2}}\OOO^\times$.  The $288$ mappings
$z\mapsto uzv^{-1}$ and $z\mapsto u\,\ov z\,v^{-1}$ are Euclidean
symmetries of $\Sigma_\OOO^\HH$, and they form the Coxeter group
$[[6,2,6]]$ of the symmetries of the $6\,$-$6$ duoprism
$\Sigma_\OOO^\HH$.

With the notation of Lemma \ref{lem:geneSLO}, the stabilizer of
$\Sigma_\OOO$ in $\SLO$ consists of the $288$ matrices
$C_{a,d}=\begin{pmatrix} a&0\\0&d\end{pmatrix}$ and
$JC_{a,d}=\begin{pmatrix} 0&d\\a&0\end{pmatrix}$ with
$a,d\in\OOO^\times$.  When $\Sigma_\OOO$ is identified with
$\Sigma_\OOO^\HH$ by the horizontal projection, the diagonal matrices
induce by Equation \eqref{eq:Poincareextension} $72$ rotational
symmetries of $\Sigma_\OOO^\HH$, and the antidiagonal ones induce
another $72$ orientation-reversing symmetries, together forming a
subgroup of index $2$ in $[[6,2,6]]]$.

It is straightforward but tedious to check that the subgroup of
$[[6,2,6]]$ that arises from the diagonal matrices in $\SLO$ acts
transitively on the vertices of $\Sigma_\OOO^\HH$.  Thus, this
subgroup acts transitively on the vertices of the fundamental cell
$\Sigma_\OOO$, which implies that all vertices of the spine $X_\OOO$
are in the same orbit.

We will now turn to a study of the link of a vertex in $X_\OOO$. Let
$$
v_0=\Big(z_0,\frac 1{\sqrt 3}\Big)\;,
$$
which is the vertex of $\Sigma_\OOO$ whose projection to $\HH$ is
$z_0$. Let $g:\HH\cup\{\infty\}\ra\HH\cup\{\infty\}$ be the homography
$z\mapsto \frac{1}{3}\,(z-z_0)^{-1}+z_0$ induced by
$$
M=\begin{pmatrix}
1 & z_0\\0 & 1
\end{pmatrix}
\begin{pmatrix}
0 & 1\\3 & -3 z_0
\end{pmatrix}
=\begin{pmatrix}
3 z_0 & 1-3z_0^2\\3 & -3 z_0
\end{pmatrix}\in\GL_2(A)\;.
$$

\blemm\label{lem:ginnormalizer}
The element $M$ belongs to the normalizer of $\SLO$ in $\SLH$.
\elemm

\dem
Computations (using Mathematica and SAGE) show that $M$ conjugates
all the generators of $\SLO$ given in Corollary \ref{cor:geneSLO} to
elements of $\SLO$, as follows. We have
$$
MJM^{-1}
%=\begin{pmatrix}\frac{7+\sqrt 3 i+3 j+\sqrt 3 k}{2}
%& -\sqrt 3i-3 j-\sqrt 3 k \\
%\frac{7-\sqrt 3 i-3 j-\sqrt 3 k}{2}
%& -\frac{7+\sqrt 3 i+3 j+\sqrt 3 k}{2}  \end{pmatrix}
=\begin{pmatrix} 3 + \rho + j + \rho j  & 1 - 2\rho - 2j - 2\rho j\\
4 - j - \rho - \rho j & - 3 - \rho - j - \rho j
\end{pmatrix}\,,
$$
$$
MT_1M^{-1}
%=\begin{pmatrix}\frac{5+\sqrt 3 i+3 j+\sqrt 3 k}{2}
%& \frac{1-\sqrt 3 i-3j-\sqrt 3 k}{2} \\
%3  &
%\frac{-1-\sqrt 3 i-3 j-\sqrt 3 k}{2}  \end{pmatrix}
=\begin{pmatrix}2 + \rho + j + \rho j & 1 - \rho -j - \rho j\\
3 &  - \rho - j - \rho j
\end{pmatrix}\,,
$$
$$
MT_jM^{-1}
%=\begin{pmatrix}\frac{-1-\sqrt 3 i+3 j+\sqrt 3 k}{2}
%& \frac{3+\sqrt 3 i-j+\sqrt 3 k}{2} \\
%3 j
%& \frac{5-\sqrt 3 i-3 j+\sqrt 3 k}{2}  \end{pmatrix}
=\begin{pmatrix} -\rho + j + \rho j & 1 + \rho -j + \rho j \\
3j &  3 - \rho - 2 j + \rho j
\end{pmatrix}\,,
$$
$$
MT_{\rho} M^{-1}
%=\begin{pmatrix}\frac{2+2 \sqrt 3 i+3 j-\sqrt 3 k}{2}
%&1-\sqrt 3 i \\
%\frac{3 +3\sqrt 3 i}2
%& 1-\sqrt 3 i-\sqrt 3 k \end{pmatrix}
=\begin{pmatrix} 2\rho + 2j - \rho j & 2 - 2\rho \\
3\rho &  2 - 2\rho +j - 2\rho j
\end{pmatrix}\,,
$$
$$
MT_{\rho j} M^{-1}
%=\begin{pmatrix}\frac{-1+\sqrt 3 i+2\sqrt 3 k}{2}
%& \frac{3+\sqrt 3 i+j-\sqrt 3 k}{2} \\
%\frac{3 j+3\sqrt 3 k}2
%& \frac{5+\sqrt 3 i-3 j-\sqrt 3 k}{2} \end{pmatrix}
=\begin{pmatrix} -1 + \rho - j + 2\rho j & 1 + \rho + j - \rho j \\
3\rho j &  2 + \rho -j - \rho j
\end{pmatrix}\,. 
$$
Since $C_{u,v}C_{u',v'}=C_{uu',vv'}$ and $JC_{u,v}J=C_{v,u}$ for all
units $u,v,u',v'$ of $\OOO$, it suffices to check the following
elements:
$$
MC_{1,-1} M^{-1}
%=\begin{pmatrix}-2+\sqrt 3 i+3 j+\sqrt 3 k
%& 4 \\
%3+\sqrt 3 i+3 j+\sqrt 3 k
%& 2-\sqrt 3 i-3 j-\sqrt 3 k  \end{pmatrix}\,,
=\begin{pmatrix} -3 + 2\rho + 2j + 2\rho j & 4 \\
2 + 2\rho + 2j + 2\rho j &  3 - 2\rho - 2j - 2\rho j
\end{pmatrix}\,,
$$
$$
MC_{1,j} M^{-1}
%=\begin{pmatrix}1+\sqrt 3 i+3 j
%& 2-\sqrt 3 i-j \\
%3+\sqrt 3 i
%& -\sqrt 3 i- j-\sqrt 3 k  \end{pmatrix}\,,
=\begin{pmatrix} 2\rho + 3j & 3 - 2\rho - j \\
2 + 2\rho &  1 - 2\rho - 2\rho j
\end{pmatrix}\,,
$$
and 
$$ 
MC_{1,\rho}M^{-1}
%=\begin{pmatrix}1+ \sqrt{3}i+\sqrt{3}k
%& 1-\sqrt{3} k \\
%\frac{3-\sqrt{3}i+2\sqrt{3}k}{2}
%& \frac{1+\sqrt 3 i-3j-\sqrt 3 k}{2} \end{pmatrix}\,.
=\begin{pmatrix} 2\rho - j + 2\rho j & 1 + j - 2\rho j \\
2 - \rho - j + 2\rho j &  \rho - j - \rho j
\end{pmatrix}\,.
$$
Thus, $M$ belongs to  the normalizer of $\SLO$ in $\SLH$.
\cqfd
%Note also that  the set $V_{3,3}$ of
%vertices of the $3\,$-$3$ duoprism %$C_{3,3}$\todo{huono merkinta} 
%lies in the sphere of
%radius $\sqrt\frac 23$ centered at $z_0$.

\bprop\label{prop:linkofv0}
If $D_A=3$, then the set of $\alpha\in A$ such that $v_0$ belongs to
the boundary of $B_\alpha(\sqrt{3})$ is
$$
V=V_{3,3}\cup g(V_{3,3})\cup\{\infty,z_0\}\;.
$$
For every $\alpha\in A$, the point $v_0$ of $\hcr$ does not belong to
the interior of $B_\alpha(\sqrt{3})$.
\eprop

The second claim implies that when $r<\sqrt{3}$, the family
$\big(B_\alpha(r)\big)_{\alpha\in A}$ does not cover $\hcr$. In
particular, the inclusions in Proposition \ref{prop:propriFVcells} (2)
are also sharp when $D_A=3$.

\medskip
\dem First observe that $v_0$ as well as all the vertices of
$\Sigma_\OOO$ are in the horizontal plane $\{(z,t)\in\hcr:t=\frac
1{\sqrt 3}\}$, which is the boundary of $B_\infty(\sqrt 3)$.

For every $\alpha\in A$, recall from Proposition \ref{prop:horoballs}
(2) that the horoball $B_{\alpha}(\sqrt 3)$ is the Euclidean ball
tangent to $\HH$ at $\alpha$ with Euclidean radius $\frac
{\sqrt{3}\;\n(I_\alpha)}{2}$. Writing $\alpha=pq^{-1}$ with
$p,q\in\OOO$ relatively prime, we have $\n(I_\alpha)=\n(q)^{-1}$. Thus
if $v_0\in B_{\alpha}(\sqrt{3})$, then the Euclidean diameter
$\sqrt{3}\n(I_\alpha)$ of $B_{\alpha}(\sqrt 3)$ is at least the
Euclidean height $\frac1 {\sqrt{3}}$ of $v_0$, that is
$\n(I_\alpha)\geq\frac{1}{3}$. Equality is only possible if $\alpha$
is the vertical projection to $\HH$ of $v_0$, that is
$\alpha=z_0$. Since $z_0=(j+\rho)(1+\rho)^{-1}$ and $j+\rho$, $1+\rho$
are relatively prime (their norms are $2$ and $3$), we have $z_0\in A$
and $\n(I_{z_0})=\frac13$. Hence the point $v_0$ does belong to the boundary
of $B_{z_0}(\sqrt 3)$, and if $\alpha\neq z_0$, then $\n(I_\alpha)=1$
or $\n(I_\alpha)=\frac{1}{2}$.

\begin{figure}[H]
\begin{center}
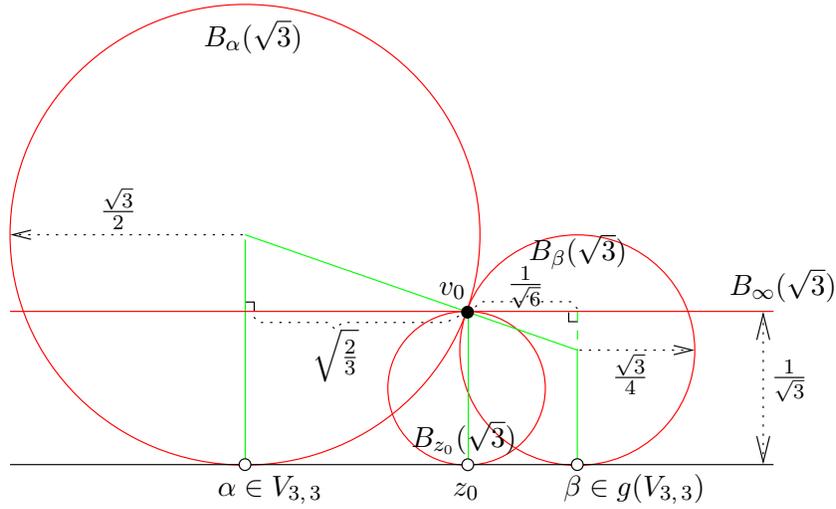
\end{center}
\caption{Intersection pattern at $v_0$ of the covering
  family of horoballs $\big(B_\alpha(\sqrt{3})\big)_{\alpha\in
    A}$}\label{fig:v0horoball}
\end{figure}

First assume that $\n(q)=1$, or equivalently that $\alpha\in\OOO$.
Then $\n(I_\alpha)=1$, hence $B_\alpha(\sqrt 3)$ is the Euclidean ball
of center $(\alpha, \frac{\sqrt 3}{2})$ and radius $\frac{\sqrt 3}
{2}$, that intersects the horizontal plane at height $\frac 1{\sqrt
  3}$ in a horizontal ball centered at $(\alpha, \frac 1{\sqrt 3})$
and of radius $\sqrt{\frac 23}$. The $9$ vertices of the fundamental
$3\,$-$3$ duoprism  of $\OOO$ are exactly at this distance
from $z_0$, and all other elements of $\OOO$ are at greater distance
from $z_0$. Hence (see Figure \ref{fig:v0horoball} on its left), $v_0$
belongs to the boundary of $B_\alpha(\sqrt 3)$ for every $\alpha\in
V_{3,3}$ and $v_0\notin B_\alpha(\sqrt 3)$ if $\alpha\in
\OOO-V_{3,3}$.

We begin the treatment of the remaining case $\n(q)=2$ by geometric
observations.  The homography $g$ defined before Lemma
\ref{lem:ginnormalizer} maps $\infty$ to $z_0$, $z_0$ to $\infty$, and
the sphere in $\HH$ of center $z_0$ and radius $r$ to the sphere in
$\HH$ of center $z_0$ and radius $\frac{1}{3r}$, for every $r>0$.  In
particular, $g$ maps the sphere in $\HH$ of center $z_0$ and radius
$\frac1{\sqrt{3}}$ to itself and the Poincar\'e extension of $g$ to
$\hcr$ (again denoted by $g$) fixes $v_0$. Thus, $g(B_\infty(\sqrt
3))=B_{z_0}(\sqrt 3)$ and Lemma \ref{lem:ginnormalizer} implies that
$g$ preserves the $\SLO$-equivariant family $(B_\alpha(\sqrt
3))_{\alpha\in A}$ of horoballs.

\medskip
Now let $\beta=pq^{-1}\in A$ be such that $v_0\in B_\beta(\sqrt{3})$
and $\n(q)=2$. Note that the Euclidean perpendicular projection from
$\hcr$ to $\HH$ does not increase the Euclidean distances, and that
the projection of the Euclidean center of $B_\beta(\sqrt{3})$ is
$\beta$ and the projection of $v_0$ is $z_0$ (see Figure
\ref{fig:v0horoball}).  Since the radius of $B_\beta(\sqrt{3})$ is
$\frac{\sqrt{3}}{4}$, we hence have $d(z_0,\beta)\leq
\frac{\sqrt{3}}{4}< \frac{1}{\sqrt{3}}$. Since $g$ fixes $v_0$ and
$gB_\beta(\sqrt{3}) =B_{g(\beta)}(\sqrt{3})$ by the above lemma, the
element $\alpha= g^{-1}(\beta)$, which satisfies $v_0\in
B_\alpha(\sqrt{3})$ and is outside the ball of center $z_0$ and radius
$\frac{1}{\sqrt{3}}$, hence cannot have a denominator of norm
$\frac{1}{2}$. Therefore $\alpha$ has denominator $1$ and by the
previous case, it belongs to $V_{3,3}$ and $v_0$ lies in the boundary
of $B_\alpha(\sqrt{3})$. So that $\beta=g(\alpha)$ belongs to
$g(V_{3,3})$ and $v_0$ lies in the boundary of $B_\beta(\sqrt{3})$.
\cqfd

\medskip
An easy computation gives
\begin{align*}
g & (V_{3,3})
=\Big\{ \frac{1+j}2 =\frac 1{1-j}\,,\; 
\frac{1+\rho j}2=\frac 1{1-j\bar\rho }\,,\; 
\frac{\rho+j}2=\frac 1{\bar\rho-j}\,,\; 
\frac{\rho(1+j)}2=\frac 1{(1-j)\bar\rho}\,,\; \\ & 
\frac{1+j+\rho j}2=\frac1{1-j(\bar\rho-1)}+j\,,\;
\frac{\rho+j+\rho j}2\,, 
\;\frac{1+\rho+j}2\,,\; 
\frac{1+\rho+\rho j}2\,,\;
\frac{1+\rho+j+\rho j}2\Big\}\;.
\end{align*}
As any element $\beta$ in $g (V_{3,3})$ is the sum of an element of
$\OOO$ with the inverse of an element of $\OOO$ with reduced norm $2$,
we have $\n(I_\beta)= \frac12$ and the horoball $B_\beta(\sqrt 3)$ has
Euclidean radius $\frac{\sqrt 3}4$.  This horoball intersects the
horizontal plane $\{(z,t)\in\hcr : t = \frac 1{\sqrt 3}\}$ in a
horizontal ball of Euclidean radius $\frac 1{\sqrt 6}$. In particular,
the points in $g (V_{3,3})$ are at Euclidean distance $\frac 1{\sqrt
  6}$ of $z_0$ and the horoballs tangent to $v_0$ are positioned as in
Figure \ref{fig:v0horoball}.

\begin{figure}[H]
\begin{center}
\includegraphics[width=12cm]{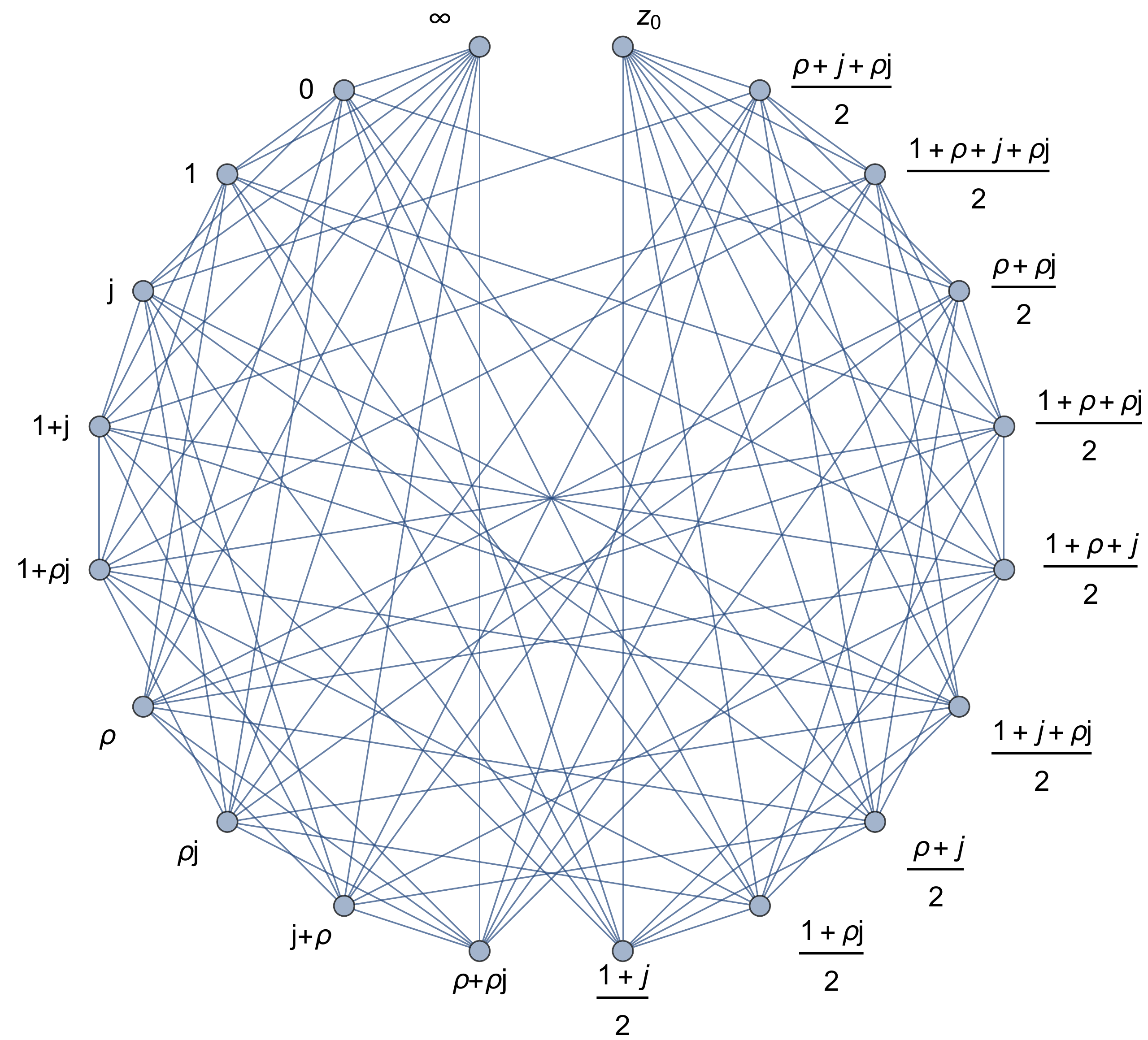}
\end{center}
\caption{Intersection pattern into
$4$-dimensional cells of Ford-Voronoi cells $\H_\alpha$ for $\alpha\in
V$.}\label{fig:Ex45intpattern}
\end{figure}

\medskip
By Proposition \ref{prop:linkofv0}, the link of $v_0$ in the
cellulation of $\hcr$ by the Ford-Voronoi cells of $\OOO$ has 20
$4$-cells, which are the intersections of a small sphere centered at
$v_0$ with the Ford-Voronoi cells $\H_\alpha$ for $\alpha$ in
$V=V_{3,3}\cup g(V_{3,3})\cup\{\infty,z_0\}$. Furthermore, for all
$\alpha\neq \beta$ in $V$, the horoballs $B_\alpha(\sqrt 3)$ and
$B_\beta(\sqrt 3)$ are tangent at $v_0$ if and only if
$\{\alpha,\beta\}$ is one of the $10$ pairs $\{\infty,z_0\}$,
$\{0,\frac{1+\rho+j+\rho j}2\}$, $\{1,\frac{\rho+j+\rho j}2\}$,
$\{\rho,\frac{1+j+\rho j}2\}$, $\{j,\frac{1+\rho+\rho j}2\}$,
$\{1+j,\frac{\rho+\rho j}2\}$, $\{1+\rho j,\frac{\rho+j}2\}$, $\{\rho
j,\frac{1+\rho+j}2\}$, $\{j+\rho,\frac{1+\rho j}2\}$ and $\{\rho+\rho
j,\frac{1+j}2\}$. By analyzing the intersections of the horoballs
$B_\alpha(1)$ contained in the Ford-Voronoi cells incident to $v_0$,
we find that each Ford-Voronoi cell containing $v_0$ intersects $9$
others in $4$-dimensional cells, that are images under $\SLO$ of the
fundamental cell $\Sigma_\OOO$, combinatorially equal to the $6\,$-$6$
duoprism $\{6\}\times\{6\}$.  The graph in Figure
\ref{fig:Ex45intpattern} shows the intersection pattern of the
$\H_\alpha$ for $\alpha\in V$.

Thus the number of ($6\,$-$6$ duoprismatic) $4$-dimensional cells of
$X_\OOO$ containing $v_0$ is exactly $90=(20\times 9)/2$, one for each
edge of this diagram.

Consider the elements $g_{\infty,1}= \begin{pmatrix} \ \,\rho &
  \rho^{-1}\\ 0 & \rho^{-1} \end{pmatrix}$, $g_{\infty,2}
= \begin{pmatrix} \ \,\rho & j\rho\\ 0 & \rho \end{pmatrix}$ and
$h_\infty= \begin{pmatrix} \;\rho j & 0 \\ 0 & j \end{pmatrix}$ in
$\SLO$ inducing respectively the homographies
$$
z\mapsto \rho z\rho+1,\;\;\;
z\mapsto \rho z\rho^{-1}+j\;\;\;{\rm and}\;\;\;
h_\infty(z)= -\rho j z j\;.
$$ 
Using the facts that $z_0=\frac{1+j+\rho+\rho j}{3}$ and $\rho j =
j\rho^{-1}$, an easy computation shows that they fix $z_0$ and
$\infty$, hence fix $v_0$ since they preserve the geodesic line
between $z_0$ and $\infty$ and the horospheres centered at $\infty$.
Hence $g_{\infty,1}$, $g_{\infty, 2}$ and $h_\infty$ belong to the
stabilizer $G_{v_0,\infty}$ of $\infty$ (or equivalently $z_0$) in the
stabilizer of $v_0$ in $\SLO$. These three elements actually generate
$G_{v_0,\infty}$.  Similar computations give that

$\bullet$~ the group $G$ generated by $g_{\infty,1}$ and $g_{\infty,
  2}$ is isomorphic to $\ZZ/3\ZZ\times \ZZ/3\ZZ$, 

$\bullet$~ $h_\infty$ has order $2$ and conjugates $g_{\infty,1}$ and
$g_{\infty, 2}$, hence each element of the abelian group $G$, to its
inverse.

\begin{figure}[H]
\begin{center}
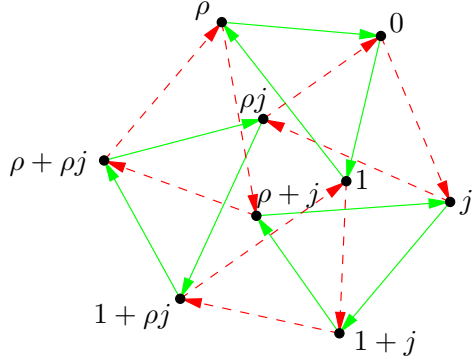
\end{center}
\caption{Transitive action of
$G_{v_0,\infty}$ on $V_{3,3}$.}\label{fig:V33action}
\end{figure}

\noindent Thus the group generated by $g_{\infty,1}$, $g_{\infty, 2}$
and $h_\infty$ is a semidirect product $(\ZZ/3\ZZ\times
\ZZ/3\ZZ)\rtimes\ZZ/2\ZZ$ with $18$ elements. The subgroup
$G_{v_0,\infty}$ acts transitively on $V_{3,3}$~: The graph in Figure
\ref{fig:V33action} shows how the points of $V_{3,3}$ are mapped by
$g_{\infty,1}$ (in continuous green) and $g_{\infty,2}$ (in dotted
red).

Since the inversion $g$ conjugates $g_{\infty,1}$ and $g_{\infty, 2}$
to their inverses, the group $G_{v_0,\infty}$ also acts transitively
on $g(V_{3,3})$.  By easy computations, the element $g_{\rho}=
\begin{pmatrix} \ \,0 & 1\\ -1 & 1 \end{pmatrix}$, inducing
the homography $z\mapsto(1-z)^{-1}$, is an element of the stabilizer
of $v_0$ in $\SLO$; it fixes $\rho\in V_{3,3}$ and $\frac{1+j+\rho j}2
\in g(V_{3,3})$, maps $\infty$ to $0\in V_{3,3}$ and $\rho j\in
V_{3,3}$ to $\frac{1+\rho j}2\in g(V_{3,3})$, and does not fix $z_0$.
Since $G_{v_0,\infty}$ acts transitively on $V_{3,3}$ and on
$g(V_{3,3})$, it follows that the stabilizer of $v_0$ acts
transitively on $V=V_{3,3}\cup g(V_{3,3})\cup \{\infty,z_0\}$.  The
stabilizer of $v_0$ in $\SLO$ coincides with the group generated by
$g_{\infty,1}$, $g_{\infty,2}$, $h_\infty$ and $g_\rho$.  It has
$20\times 18=360$ elements (the number of $4$-cells of the link of
$v_0$ in the tesselation of $\hcr$ times the order of the stabilizer
of one $4$-cell, the one corresponding to $\H_\infty$).

The dual tiling of the $6$-$6$-duoprismatic tiling of $\HH$ is the
$3$-$3$-duoprismatic tiling. Therefore, the link of $v_0$ in $\partial
\H_\infty$ (hence in all $\partial \H_\alpha$ containing $v_0$) is the
$3$-skeleton of the dual of the $3$-$3$ duoprism, namely the $3$-$3$
duopyramid, whose Schl\"afli symbol is $\{3\}+\{3\}$ and whose
symmetry group has order $8\times 3^3=72$.
% that has $9$ tetrahedral $3$-cells, $18$ triangles, $15$ edges
%and $6$ vertices.
The group generated by $g_{\infty,1}$, $g_{\infty,2}$ and $h_\infty$
is a subgroup of index $4$ in the full group of symmetries of the link
of $v_0$ in $\partial \H_\infty$.  The link of $v_0$ in $\hcr$ is
constructed of $20$ copies of
%the  above complex of $9$ tetrahedra 
the $3$-$3$ duopyramid, that are glued together according to the
intersection pattern described above, forming the $4$-skeleton of the
dual of the {\em birectified $5$-simplex}. Since the birectified
$5$-simplex is called the {\em dodecateron} and has twelve $4$-faces,
its dual, which has twenty $4$-faces and does not seem to have a name
in the literature, could be called the {\em icosateron}.  The full
group of symmetries %(including orientation-reversing ones) 
of the dual of the icosateron, whose Coxeter notation is $[[4^3]]$,
has $1440=4\times 360$ elements. The stabilizer of $v_0$ in $\SLO$ is
naturally identified with a subgroup of index $4$ in $[[4^3]]$.  This
concludes the study of Example \ref{ex:eisenspine}.

\begin{comment}
  Thus, the link has $\frac{20\cdot 9}2 =90$ tetrahedra.  To make
  things concrete, let us concentrate on the tetrahedron
  $T_{\infty,0}=\partial\H_\infty\cap\partial\H_0$. The boundary of
  $T_{\infty,0}$ consists of $4$ triangles,
  $$
  \Delta_{\infty,0,1}=\partial\H_\infty\cap\partial\H_0\cap \partial\H_1
  =T_{\infty,0}\cap T_{0,1}\cap T_{\infty,1}\,,
  $$ and the similarly defined $\Delta_{\infty,0,\rho}$,
  $\Delta_{\infty,0,j}$, $\Delta_{\infty,0,\rho j}$. Each triangle in
  the link of $v_0$ in $X_\OOO$ is contained in $3$ tetrahedra and,
  therefore, there are $\frac{20\cdot 18}3 =120$ triangles in the
  link.  \todo{\tiny I guess that the link is the $3$-skeleton of the
    dual of the dodecateron with $60$ edges and $12$ vertices but I am
    not able to prove it for the moment.}
\end{comment}

\eexem

\rem When $D_A\in\{2,3,5\}$, let $P_\OOO$ be the hyperbolic
$5$-polytope that consists of the points in the halfspace $\H_{\infty
  0}$ whose horizontal projection to $\HH$ is $\Sigma_\OOO^\HH$. The
quotient orbifold $\SLO\bs\hcr$ is obtained from $P_\OOO$ by gluing
the vertical sides of $P_\OOO$ by the translations in the stabilizer
of $\infty$ in $\SLO$, and then folding by the action of the
stabilizer of $\Sigma_\OOO$. The quotient space $\SLO\bs X_\OOO$
obtained by making the above identifications in $\Sigma_\OOO$ is a
$4$-dimensional cellular retract of $\SLO\bs\hcr$ that could be used to
study the homology of $\SLO$ and $\PSL_2(\OOO)$ analogously to the
study of the Bianchi groups in \cite{Mendoza80} and \cite{SchVog83}.

%\todo{\tiny to be completed for $D_A=5$ hence
%  $A=\big(\frac{-2,-5}{\QQ}\big)$
% $\OOO=\ZZ[1, (1 + i + j)/2, j, (2 + i + k)/4]$ \cite{Vigneras80} p. 98  ?}
%% 

\section{Waterworlds}
\label{sec:waterworld}

Let $A$ be a definite quaternion algebra over $\QQ$ and let $\OOO$ be
a maximal order in $A$. Let $f$ be an indefinite integral binary
Hamiltonian form over $\OOO$.

The form $f$ defines a function $F=F_f:\PP^1_r(A)\ra\QQ$ by
$$
F([x:y])=\frac{f(x,y)}{\n(\OOO x+\OOO y)}\;.
$$ 
This definition does not depend on the choice of representatives
$(x,y) \in A\times A$ of $[x:y]\in\PP^1_r(A)$, and $f$ is uniquely
determined by its associated function $F$. In particular, we may take
$x,y\in\OOO$ in order to compute $F([x:y])$, so that the numerator of
the fraction defining $F([x:y])$ belongs to $\ZZ$. Note that
$\SL(\OOO)$ acts with finitely many orbits on $\PP^1_r(A)$, since the
number of cusps is finite, and that the denominator defining
$F([x:y])$ is invariant under $\SL(\OOO)$. Therefore there exists
$N\in\NN-\{0\}$ such that $F$ has values in $\frac{1}{N}\ZZ$, hence
the set of values of $F$ is discrete.

Note that for every $g\in\SLO$, the function $F_{f\circ g}$ associated
to the form $f\circ g$ is $F\circ g$ (where we again denote by $g$ the
projective transformation of $\PP^1_r(A)$ induced by $g$). In
particular, $F\circ g=F$ if $g\in\automH$.

\begin{figure}[H]
\begin{center}
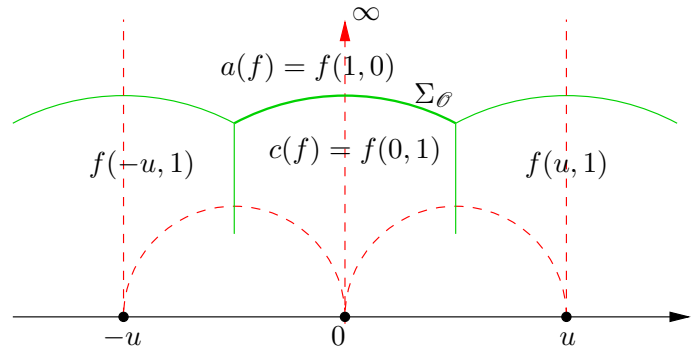
\end{center}
\caption{Values of $F$ on
Ford-Voronoi cells meeting $\Sigma_\OOO$.}\label{fig:valuesofF}
\end{figure}

As in \cite{Conway97} for integral indefinite binary quadratic forms,
we will think of $F$ as a map which associates a rational number to
(the interior of) any Ford-Voronoi cell. For instance, if $D_A=2$ and
$\OOO$ is the Hurwitz order, then the values of $F$ on the two
Ford-Voronoi cells $\H_\infty,\H_0$ containing the fundamental cell
$\Sigma_\OOO$ are $f(1,0), f(0,1)$ and the values of $F$ on the $24$
Ford-Voronoi cells meeting $\Sigma_\OOO$ in a $3$-dimensional cell are
$f(u,1)$ for $u\in\OOO^\times$ (see Figure \ref{fig:valuesofF}).
%\addtocounter{fig}{1}\arabic{fig}\addtocounter{fig}{-1}).

%\begin{center}
%\input{fig_valuesFonFVcell.pdf_t}\\
%\addtocounter{fig}{1} Figure \arabic{fig}: Values of $F$ on
%Ford-Voronoi cells meeting $\Sigma_\OOO$.
%\end{center}

Let $\mmm$ be a left fractional ideal of $\OOO$.  For every $s\geq 0$,
let
$$
\psi_{F,\,\mmm}(s)=\card\;\;_{\mbox{$\automH$}}\bs
\big\{(u,v)\in\mmm\times\mmm\;:\;|F(u,v)|\leq s,
\;\;\;\OOO u+\OOO v=\mmm\big\}\;,
$$ 
which is the number of nonequivalent $\mmm$-primitive representations
by $F$ of rational numbers in $\frac 1N\ZZ$ with absolute value at
most $s$. We showed in \cite[Theo.~1]{ParPau13ANT} and
\cite[Cor.~5.6]{ParPau14AFST} that there exists $\kappa>0$ such that,
as $s$ tends to $+\infty$,
$$
\psi_{f,\,\mmm}(s)= \frac{45\;D_A\;\covol(\automH)}
{2\,\pi^2\;\zeta(3)\;\Delta(f)^2\;\prod_{p|D_A}(p^3-1)}\;\;
s^4(1+\operatorname{O}(s^{-\kappa}))\;.
$$

\blemm \label{lem:existflooded}
The function $F$ takes all signs $0, +, -$. 
\elemm

\dem It takes positive and negative values since $f$ is
indefinite. The values of $F$ are actually positive at the points in
$\PP^1_r(A)$ in one of the two components of $\PP^1_r(\HH)-
\C_\infty(f)$ and negative at the ones in the other component. But
contrarily to the cases of integral binary quadratic and Hermitian
forms, all integral binary Hamiltonian forms $f$ over $\OOO$ represent
$0$, since by taking a $\ZZ$-basis of $\OOO$, the form $f$ becomes an
integral binary quadratic form over $\ZZ$ with $8$ variables and all
integral binary quadratic forms over $\ZZ$ with at least $5$ variables
represent $0$, see for instance \cite[p.~77]{Serre70} or
\cite[p.~75]{Cassels78}.  \cqfd

\medskip
A Ford-Voronoi cell will be called {\it flooded} for $f$ if the value
of $F$ on its point at infinity is $0$. Lemma \ref{lem:existflooded}
says that there are always flooded Ford-Voronoi cells. See also
\cite[Cor.~4.8]{Vulakh93}. The flooded Ford-Voronoi cells for $f$
correspond to Conway's {\it lakes} for an isotropic integral
indefinite binary quadratic form over $\ZZ$, see \cite[page
  23]{Conway97}. There were only two lakes, whereas there are now
countably infinitely many flooded Ford-Voronoi cells for $f$, one for
each parabolic fixed point of the group of automorphs of $f$.

\bexem 
Consider the definite quaternion algebra $A$ with $D_A=2$, $\OOO$ the
Hurwitz order and a Hamiltonian form $f$ with $a(f)=0$, $b=b(f),c=c(f)
\in \ZZ-\{0\}$ such that $b$ does not divide $c$ nor $2c$.  Then
$\H_\infty$ is flooded. Let $\alpha=xy^{-1}$ with $x\in\OOO$ and
$y\in\OOO-\{0\}$ relatively prime. If $\n(y)\le 2$, then the
Ford-Voronoi cell $\H_\alpha$ is not flooded, since otherwise the
equation $b\tr(\bar x\,y)+c\n(y)=0$ would imply that $b$ divides $c$
or $2c$.  If $\n(y)>2$, then $\n(I_\alpha)=\frac{\n(\OOO x+\OOO
  y)}{\n(y)}=\frac{1}{\n(y)}<\frac12$. Hence by Proposition
\ref{prop:horoballs} (2), we have $B_\alpha(\sqrt{2})\cap
B_\infty(\sqrt{2})=\emptyset$. Therefore
$\H_\alpha\cap\H_\infty=\emptyset$ by Proposition
\ref{prop:propriFVcells} (2). This proves that $\H_\infty$ does not
meet any other flooded Ford-Voronoi cell. Thus if
%%
%\todo{when is this true ?}
%%
the hyperbolic $4$-orbifold $\automH\bs \C(f)$ has only one cusp, then
the flooded Ford-Voronoi cells are pairwise disjoint.  We actually do
not know when $\automH\bs \C(f)$ has only one cusp. \eexem

We have the following analog of  the statement of Conway (loc.~cit.)
that the values of the binary quadratic form along a lake are in an
infinite arithmetic progression.

\bprop \label{prop:arithprogress} Let $\alpha_0\in\PP^1_r(A)$ be such
that the Ford-Voronoi cell $\H_{\alpha_0}$ is flooded for $f$.  If
$\alpha_0$ belongs to the $\SLO$-orbit of $\infty$, let
$\Lambda_{\alpha_0}= \OOO$.  Otherwise, let
$$
\Lambda_{\alpha_0}=\OOO\cap \alpha_0^{-1}\OOO \cap
\OOO\alpha_0^{-1}\cap \alpha_0^{-1}\OOO\alpha_0^{-1}\;.
$$ 
Then there exists a finite set of nonconstant affine maps
$\{\varphi_j:\HH\ra \RR\;:\;j\in J'\}$ defined over $\QQ$ such that the
set of values of $F$ on the Ford-Voronoi cells meeting $\H_{\alpha_0}$
is  $\bigcup_{j\in J'}\varphi_j(\Lambda_{\alpha_0})$.  
\eprop

\dem For every $\alpha\in\PP^1_r(A)$, let $E_\alpha=\{\beta\in
\PP^1_r(A) -\{\alpha\} \;:\; \H_\alpha\cap \H_\beta\neq \emptyset\}$.
Note that $E_{g\cdot\alpha} = g\cdot E_\alpha$ for every $g\in\SLO$,
by Proposition \ref{prop:propriFVcells} (1). 

First assume that $\alpha_0$ belongs to the $\SLO$-orbit of $\infty$.
Then up to replacing $f$ by $f\circ g$ for some $g\in\SLO$ such that
$g\cdot \infty = \alpha_0$, we may hence assume that $\alpha_0=
\infty$.

Let $a=a(f)$, $b=b(f)$ and $c=c(f)$. Note that $\H_\infty$ is flooded
for $f$ if and only if $f(0,1)=0$, that is, if and only if $a=0$. We
then have $b\neq 0$ since $f$ is indefinite. Hence $F(E_\infty)=
\big\{\frac{\tr(\overline{u}\,b)+c}{\n(I_u)}\;:\; u\in E_\infty
\big\}$.  Since the stabilizer of $\infty$ in $\SLO$ acts with
finitely many orbits on the cells of $\partial \H_\infty$, its finite
index subgroup $\OOO$ acts by translations with finitely many orbits
on $E_\infty$. Hence there exists a finite subset $J'$ of $A$ such that
$E_\infty =J'+\OOO$.  Since $I_{\alpha+o}=I_\alpha$ for all $\alpha\in
A$ and $o\in\OOO$, the result follows with $\varphi_j : u \mapsto
\frac{\tr(\overline{b}(j+u))+c} {\n(I_j)}$ for all $j\in J'$.

Assume now that $\alpha_0$ does not belong to the $\SLO$-orbit of
$\infty$, so that in particular $\alpha_0\in A-\{0\}$.  Let
$\Ga_{\alpha_0}$ be the stabilizer of $\alpha_0$ in $\SLO$, which acts
with finitely many orbits on $E_{\alpha_0}$. Let $g=
\Big(\begin{array}{cc} \alpha_0 & -1 \\ 1 & 0 \end{array}\Big)$, which
belongs to $\SL_2(A)$ and whose inverse projectively maps $\alpha_0$
to $\infty$.  Then (see for instance \cite[\S 5]{ParPau13ANT}),
$\Lambda_{\alpha_0}$ is a $\ZZ$-lattice in $\HH$, such that the group
of unipotent upper triangular matrices with coefficient $1$-$2$ in
$\Lambda_{\alpha_0}$ is a finite index subgroup of $g^{-1}
\Ga_{\alpha_0} g$.  A similar argument concludes.
%(except that $f\circ g$ has rational instead of integral coefficients).
%%
%\todo{a bit too fast}
%% 
\cqfd

\medskip 
By a {\it projective real hyperplane} in $\partial_\infty\HH^5_\RR=
\PP^1_r(\HH)=\HH\cup\{\infty\}$, we mean in what follows the boundary
at infinity of a hyperbolic hyperplane in $\HH^5_\RR$. The ones
containing $\infty= [1:0]$ are the union of $\{\infty\}$ with the
affine real hyperplanes in $\HH$. The ones not containing $\infty$ are
the Euclidean spheres in the affine Euclidean space $\HH$.

\blemm \label{lem:defuniqnonclineated}
The form $f$ is uniquely determined by the values of its
associated function $F$ at six points in $\PP^1_r(A)$ that do not lie
in a projective real hyperplane.  
\elemm

\dem Let $a=a(f)$, $b=b(f)$ and $c=c(f)$. Let first prove that we may
assume that the six points in $A\cup \{\infty\}$ are $\infty=[1:0]$,
$0$, $\alpha_0=1$ and $\alpha_1, \alpha_2, \alpha_3\in A-\{0\}$.

Note that for all $x,y\in A$ and $g\in \GL_2(A)$, if $g_1,g_2$ are the
components of the linear selfmap $g$ of $A\times A$, then
\begin{equation}\label{eq:equivF}
F_{f\circ g}([x:y])= F_f\circ g([x:y])\;
\frac{\n(\OOO g_1(x,y)+\OOO g_2(x,y))}{\n(\OOO x +\OOO y)}\;.
\end{equation}
Given six points in $\PP^1_r(A)$ not in a projective real hyperplane
of $\PP^1_r(\HH)$, the first three of them constitute a projective
frame of the projective line $\PP^1_r(A)$. Hence by the existence part
of the fundamental theorem of projective geometry (see
\cite[Prop.~4.5.10]{Berger77}), there exists an element $g\in
\GL_2(A)$ mapping them to $\infty,0,1$. Note that this existence part
does hold in the noncommutative setting, though the uniqueness part
does not. The initial claim follows by Equation \eqref{eq:equivF}.

Now, the values of $F$ at the points $\infty$, $0$, $\alpha_0,
\alpha_1, \alpha_2, \alpha_3$ give a system of six equations on the
unknown $a,b,c$, of the form $a=A_1$, $c=A_2$, $a+\tr b+c=A_3$,
$\tr(\overline{\alpha_i}\,b) =A_{i+3}$ for $i\in\{1,2,3\}$. Thus $a$
and $c$ are uniquely determined, and $b$ belongs to the intersection
of four affine real hyperplanes in $\HH$ orthogonal to
$\alpha_0,\alpha_1,\alpha_2, \alpha_3$ with equations
$\tr(\overline{\alpha_i}\,b)=A'_{i}$ for $i\in\{0,1,2,3\}$. The result
follows since if $\alpha_0,\alpha_1, \alpha_2,\alpha_3$ are linearly
independent over $\RR$, then for all $A'_0,A'_1,A'_2,A'_3\in\RR$, such
an intersection contains one and only one point of $\HH$.  \cqfd

\bprop\label{prop:funiquedeterm} Let $v$ be a vertex of the spine
$X_\OOO$. The form $f$ is uniquely determined by the values of its
associated function $F$ on the Ford-Voronoi cells containing
$v$, that is, on the points $\alpha\in \PP^1_r(A)$ such that
  $v\in\H_\alpha$.
\eprop

\dem A dimension count shows that there are at least six Ford-Voronoi
cells meeting at each vertex $v$ of the spine. Their points at
infinity cannot all be on the same projective real hyperplane
$P$. Otherwise, the intersection of the equidistant hyperbolic
hyperplanes between the pair of them yielding a $4$-dimensional cell
containing $v$ would have dimension at least $1$: It would contain a
germ of the orthogonal through $v$ to the convex hull of $P$ in
$\hcr$. The result follows by Lemma \ref{lem:defuniqnonclineated}.
\cqfd

\medskip
The {\em waterworld} of $f$ is
$$
\W(f)=\bigcup_{\alpha\ne\beta\in\PP_r^1(A),\ F(\alpha)F(\beta)<1}\H_\alpha\cap\H_\beta\,.
$$ 
As  $F\circ g=F$ if $g\in\automH$, the waterworld $\W(f)$ is  invariant under the
group of automorphs $\automH$ of $f$.

Since $f$ is always isotropic over $A$, the arguments of Conway and
Bestvina-Savin for the anisotropic case no longer apply, and the
waterworld of $f$ could be empty. We do not know precisely when the
waterworlds are nonempty, and we now study some examples.

\bexem The binary Hamiltonian form $f(u,v)=\tr(\ov u\, v)$ is
indefinite with discriminant $1$. The coefficients of $f$ are rational
integers so it is integral over any maximal order $\OOO$ of any
definite quaternion algebra $A$ over $\QQ$.  Let us prove that the
waterworld $\W(f)$ is not empty.

It is easy to check that $\C_\infty(f)=\{z\in\HH:\tr z=0\}\cup
\{\infty\}$.  Let $a\in \OOO$ be such that $\tr(a)=1$ (which does
exist since $\OOO$ is maximal, hence $\tr:\OOO\ra\ZZ$ is onto, see for instance the proof of Proposition 16 in \cite{ChePau19}). In
particular $a\neq 0$, $a\neq -\bar a$, and $a,-\bar a$ are in two
different components of $\partial_\infty\hcr-\C_\infty(f)$, so that
$F(a)F(-\bar a)<0$. Let us prove that $\H_a$ and $\H_{-\bar a}$
intersect in a $4$-dimensional cell of $X_\OOO$, which thus belongs to
$\W(f)$. By Proposition \ref{prop:propriFVcells} (2), it is sufficient
to prove that $B_a(1)$ and $B_{-\bar a}(1)$ meet. By Theorem
\ref{theo:coverhoroball} (1), this is equivalent to proving that $I_a
I_{\bar a} = \OOO (\tr a)$. But this holds since $\tr a =1$ and
$I_b=\OOO$ when $b\in \OOO$.

%Assume that $\OOO$ is left-Euclidean. In particular, $A$ is invariant
%under the mapping $z\mapsto-\overline z$ that fixes $\C_\infty(f)$
%pointwise.  Let $\alpha\in A$ with $\tr\alpha\ne 0$.  By symmetry,
%$\C(f)$ is the equidistant hyperplane of $B_\alpha(1)$ and
%$B_{-\overline\alpha}(1)$.  The Ford-Voronoi cells $\H_\alpha$ and
%$\H_{-\overline\alpha}$ intersect if and only if $B_\alpha(1)$ and
%$B_{-\overline\alpha}(1)$ are tangent.
% 
%As $\OOO$ is left-Euclidean, the group $\SLO$ acts transitively on $A$
%and on the collection of horoballs $(B_\alpha(1))_{\alpha\in
%  A\cup\{\infty\}}$. Thus, for every $\alpha\in A$ there exists
%$g=\begin{pmatrix} a & b\\c & d\end{pmatrix}\in\SLO$ such that
%$\alpha=g\cdot\infty=ac^{-1}$.  As $g\cdot(-c^{-1}d)=\infty$,
%$g\cdot(-c^{-1}d,1)=(\alpha,\frac 1{\n(c)})$ is the highest point in
%$B_\alpha(1)$. Thus the Euclidean radius of $B_\alpha(1)$ is $\frac
%1{2\n(c)}$. As the Euclidean distance of $\alpha$ from $\C_\infty(f)$
%is $|\frac{\tr\alpha}2|$, this implies that $B_\alpha(1)$ intersects
%$\C(f)$ if and only if $\tr ac=\pm 1$ which is satisfied for example
%if $\alpha$ is an integer with trace $\pm 1$.  All left-Euclidean
%maximal orders contain such elements, and thus, the waterworld of $f$
%is not empty in any of these cases.

\begin{figure}[H]
\begin{center}
\includegraphics[width=13cm]{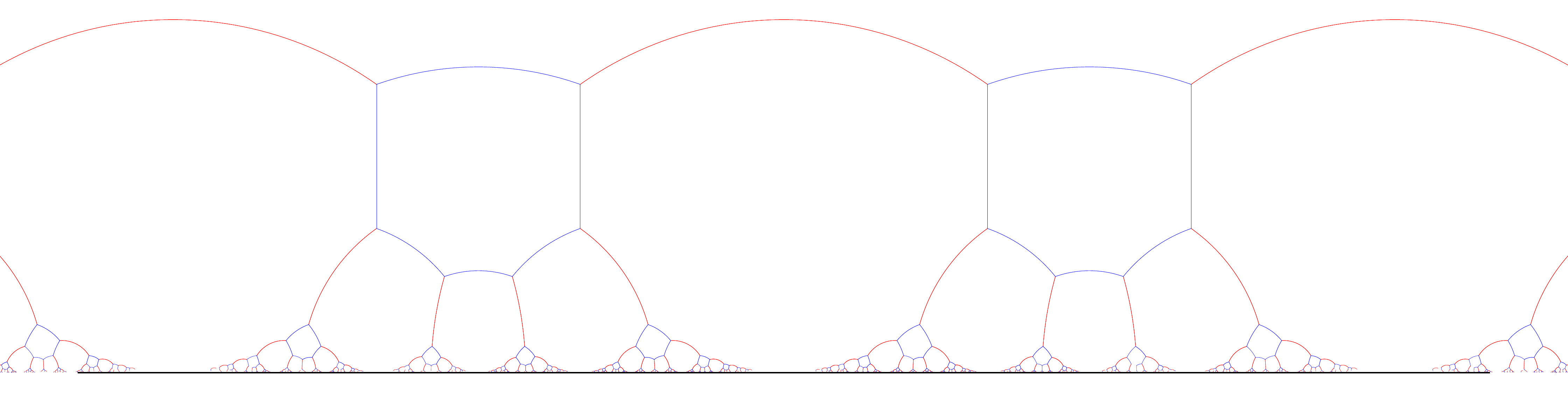}
\end{center}
\caption{Ocean
of Hermitian form $f(u,v)=\tr(\ov u \,v)$ over $\ZZ[\frac{1+i\sqrt
    3}2]$.}\label{fig:Eisensteinocean}
\end{figure}

\medskip
Figure \ref{fig:Eisensteinocean}
%\addtocounter{fig}{1}\arabic{fig}\addtocounter{fig}{-1}
illustrates the analogous case of the ocean in $\htr$ of the isotropic
binary Hermitian form $f(u,v)=\tr(\ov u \,v)$ considered as an
integral form over the Eisenstein integers $\ZZ[\frac{1+i\sqrt 3}2]$.
The blue hexagons are the components of the ocean of $f$ in the
hyperplane $\C(f)=\{(z,t)\in\htr:\Im \;z=0\}$ which is a copy of the
(upper halfplane model of the) real hyperbolic plane.

%\begin{center}
%\includegraphics[width=13cm]{fig_EWW.pdf}\\
%\addtocounter{fig}{1} Figure \arabic{fig}: Ocean
%of Hermitian form $f(u,v)=\tr(\ov u \,v)$ over $\ZZ[\frac{1+i\sqrt
%    3}2]$.
%\end{center}

\begin{figure}[H]
\begin{center}
\includegraphics[width=13cm]{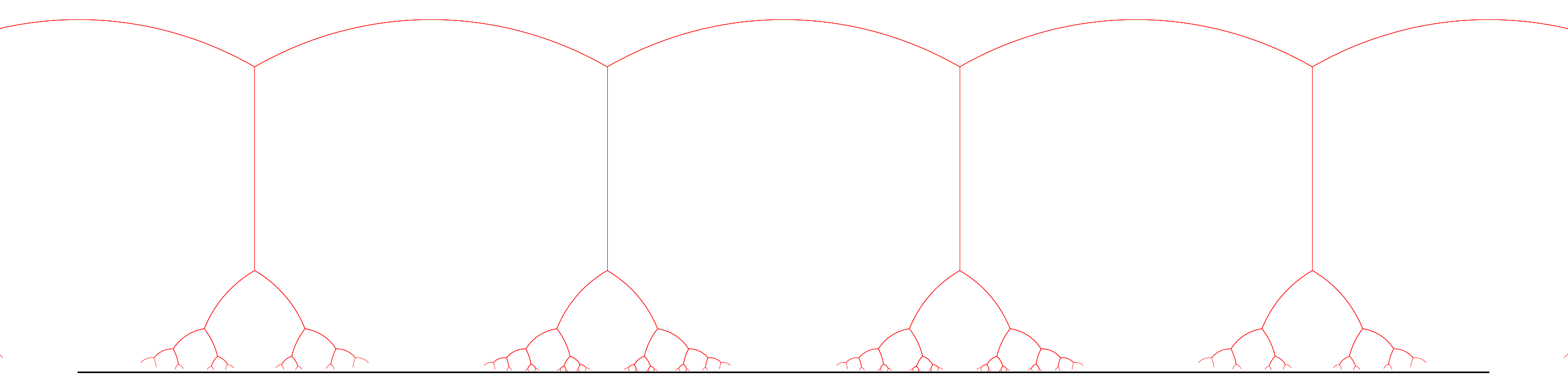}
\end{center}
\caption{Empty ocean of Hermitian
form $f(u,v)=\tr(\ov u \,v)$ over $\ZZ[i]$.}\label{fig:Gaussocean}
\end{figure}

We do not have an example of an empty waterworld and, in fact, it may
be that no such example exists. However, the ocean of the isotropic
binary Hamiltonian form $f(u,v)=\tr(\ov u \,v)$ considered over the
Gaussian integers $\ZZ[i]$ is empty (see Figure \ref{fig:Gaussocean}). 
%\addtocounter{fig}{1}\arabic{fig}\addtocounter{fig}{-1}). 
In order
to prove this, let $\alpha\in \QQ(i)$ with $\tr\alpha\ne 0$. Note that
in the commutative case, $\n(I_\alpha) = \n(I_{-\overline{\alpha}})$,
so that the Euclidean balls $B_\alpha(1)$ and
$B_{-\overline\alpha}(1)$ have the same radius.  By symmetry, $\C(f)$
is the equidistant hyperbolic hyperplane of $B_\alpha(1)$ and
$B_{-\overline\alpha}(1)$.  Since $\ZZ[i]$ is Euclidean, the spine of
$\SL_2(\ZZ[i])$ has only one orbit of $2$-cells (see
\cite{BesSav12}). Hence all the intersections of the Ford-Voronoi cells 
are in the orbit of the fundamental cell, and therefore, 
$\H_\alpha$ and
$\H_{-\overline\alpha}$ intersect if and only if $B_\alpha(1)$ and
$B_{-\overline\alpha}(1)$ are tangent, that is, if and only if
$B_\alpha(1)$ intersects $\C(f)$.

Since the hyperbolic $3$-orbifold $\SL_2(\ZZ[i])\bs\htr$ has only one
cusp, there exists $g=\begin{pmatrix} a & b\\c & d\end{pmatrix} \in
\SL_2(\ZZ[i])$ such that $\alpha= g\cdot\infty=ac^{-1}$.  Since
$g\cdot(-c^{-1}d)=\infty$, the point $g\cdot(-c^{-1}d,1)=(\alpha,\frac
1{\n(c)})$ is the highest point in $B_\alpha(1)=gB_\infty(1)$. Thus
the Euclidean radius of $B_\alpha(1)$ is $\frac 1{2\n(c)}$. As the
Euclidean distance of $\alpha$ from $\C_\infty(f)$ is
$|\frac{\tr\alpha}2|$, this implies that $B_\alpha(1)$ intersects
$\C(f)$ if and only if $\big|\frac{\tr\alpha}{2}\big|\leq
\frac{1}{2\n(c)}$, that is, if and only if $\tr a\,\bar c=\pm 1$. This
is impossible since the trace of any Gaussian integer is even.

%\begin{center}
%\includegraphics[width=13cm]{fig_GaussEmptyOcean.pdf}\\ 
%\addtocounter{fig}{1} Figure \arabic{fig}: Empty ocean of Hermitian
%form $f(u,v)=\tr(\ov u \,v)$ over $\ZZ[i]$.
%\end{center}

\eexem

\bprop If the union of the flooded Ford-Voronoi cells does not
separate $\hcr$, and in particular if the flooded Ford-Voronoi cells
are pairwise disjoint, then the waterworld of $f$ is nonempty.  
\eprop

\dem The assumption says that the topological space 
$$
X=\hcr- \bigcup_{\alpha\in \PP_r^1(A),\;F(\alpha)=0}\H_\alpha
$$ 
is connected. If $\W(f)= \emptyset$, then 
$$
X=\Big(\bigcup_{\alpha\in
  \PP_r^1(A), \;F(\alpha)<0} \H_\alpha\Big) \cup
\Big(\bigcup_{\alpha\in \PP_r^1(A),\;F(\alpha)>0} \H_\alpha \Big)
$$
would be a partition into two nonempty (since $f$ is indefinite)
locally finite, hence closed, unions of closed polyhedra,
contradicting the connectedness of $X$.  \cqfd

\bprop \label{prop:waterworldmodautomcompact}
The quotient $\automH\bs\W(f)$ is compact, and the set of
flooded Ford-Voronoi cells consists of finitely many $\automH$-orbits.
\eprop

\dem The points at infinity of the flooded Ford-Voronoi cells are the
parabolic fixed points of $\SLO$ contained in $\C_\infty(f)$, hence
are the parabolic fixed points of the group of automorphs $\automH$.
Since $\automH$ is a lattice in the real hyperbolic $4$-space $\C(f)$,
the quotient $\automH\bs\C(f)$ has only finitely many cusps. This
proves the second claim.

\medskip
Let $\alpha,\beta\in\PP_r^1(A)$ be such that $F(\alpha)F(\beta)<0$ and
the intersection $\H_\alpha \cap\H_\beta$ is nonempty. Then the
intersection $B_\alpha(\sqrt{D_A})\cap B_\beta(\sqrt{D_A})$ is
nonempty by Proposition \ref{prop:propriFVcells} (2), hence the
hyperbolic distance between the horoballs $B_\alpha(1)$ and
$B_\beta(1)$ is at most $\ln D_A$. By Lemma
\ref{lem:distBonealphabeta}, we hence have
$\frac{\n(\alpha-\beta)}{\n(I_\alpha I_\beta)}\leq D_A$.

Let $a=a(f)$, $b=b(f)$, $c=c(f)$ and $\Delta=\Delta(f)$. Write
$\alpha=[x:y]$ and $\beta=[u:v]$ with $x,y,u,v\in\OOO$ and
$y,v\in\ZZ$. Note that
$$ 
\begin{pmatrix} x & u\\ y & v\end{pmatrix}^*
\begin{pmatrix} a & b\\ \ov{b} & c\end{pmatrix}
\begin{pmatrix} x & u\\ y & v\end{pmatrix}=
\begin{pmatrix} f(x,y) & z\\ \ov{z} & f(u,v)\end{pmatrix}\;,
$$ 
for some $z\in\OOO$.  Since $y,v\in\RR$, an easy computation of
Dieudonné determinants thus gives
$$
\big|\n(z)-f(x,y)f(u,v)\big|=\n(xv-uy)\,\Delta\;.
$$ 
Hence $0\leq -f(x,y)f(u,v)\leq \n(z)-f(x,y)f(u,v)=\n(xv-uy)\,\Delta$
and
$$
0\leq -F(\alpha)F(\beta)=
\frac{-f(x,y)f(u,v)}{\n(\OOO x+\OOO y)\n(\OOO u+\OOO v)}\leq  
\frac{\n(\alpha-\beta)}
{\n(I_\alpha)\n(I_\beta)}\,\Delta\leq D_A\,\Delta\;.
$$ 
Since the set of values of $F$ is discrete in $\RR$, this implies
that $F$ takes only finitely many values on the Ford-Voronoi cells
that intersect $\W(f)$.

Given any vertex $v \in \W(f)$, for every $g\in\SLO$, if $F(\alpha)=
F(g\cdot\alpha)$ for all $\alpha\in A$ such that the Ford-Voronoi cell
$\H_\alpha$ contains $v$, then $f=f\circ g$ by Proposition
\ref{prop:funiquedeterm}. Since there are only finitely many orbits of
$\SLO$ on the vertices of the spine $X_\OOO$ and since $F$ takes only
finitely many values on the Ford-Voronoi cells meeting the waterworld
$\W(f)$, this implies that $\automH$ has only finitely many orbits of
vertices in $\W(f)$. The result follows.  
\cqfd

\medskip
Note that there exist a positive constant and finitely many pairs
$\{\alpha,\beta\}$ in $A$ such that, for all indefinite integral
binary Hamiltonian forms $f$ over $\OOO$ up to the action of $\SLO$,
the distance between $\H_\alpha$ and $\H_\beta$ is at most this
constant and $F(\alpha) F(\beta)<0$. This follows, even if the
waterworld $\W(f)$ could be empty, from the fact that the flooded
Ford-Voronoi cells only have their points at infinity on the
$3$-sphere $\C_\infty(f)$ in $\PP^1_r(\HH)$, and by the cocompactness
of the action of $\SLO$ on its spine $X_\OOO$. The above arguments
hence allow to give another proof of Corollary 25 in
\cite{ParPau13ANT}, saying that the number of $\SLO$-orbits in the set
of indefinite integral binary Hamiltonian forms over $\OOO$ with given
discriminant is finite.

\medskip
We introduce two variants of $\W(f)$. The {\em sourced waterworld} 
$\W_+(f)$ of $f$ is the union of its
waterworld and of its flooded Ford-Voronoi cells
$$
\W_+(f)=\W(f)\;\cup\bigcup_{\alpha\in\PP_r^1(A),\ F(\alpha)=0}\H_\alpha\,.
$$ 
The {\em coned-off waterworld} $\C\W(f)$ of $f$ is obtained from
$\W(f)$ by adding geodesic rays from its boundary points to the points
at infinity of the corresponding flooded Ford-Voronoi cells
$$
\C\W(f)=\W(f)\;\cup
\bigcup_{\alpha\in\PP_r^1(A),\,x\in\W(f)\cap\H_\alpha\;:\;F(\alpha)=0}
[x,\alpha[\,.
$$ 
Both the sourced waterworld $\W_+(f)$ and
the coned-off waterworld $\C\W(f)$ of $f$ are invariant under the
group of automorphs $\automH$ of $f$.

\bigskip Before stating the main result of this paper, we give two
lemmas and refer to Section 6 of \cite{BesSav12}  for the proofs

\blemm \label{lem:propproj} Let $P,P'$ be hyperbolic hyperplanes in
$\HH^n_\RR$ that do not intersect perpendicularly. Then the closest
point mapping from $P$ to $P'$ is a homeomorphism onto a convex open
subset of $P'$, which maps any hyperbolic polyhedron of $P$ to a
hyperbolic polyhedron of $P'$. \cqfd 
\elemm

\blemm \label{lem:claim3} Let $f$ be an indefinite integral binary
Hamiltonian form over $\OOO$. If $\ell$ is a geodesic line in $\hcr$ that is 
perpendicular to the halfplane $\C(f)$, oriented such that $\ell(\pm\infty) \in
\{[x:y]\in\PP^1_r(\HH)\;:\; \pm f(x,y)>0\}$, if $\ell$ meets
transversally at a point $z$ the interior of a $4$-dimensional cell
$\H_{\alpha_-}\cap \H_{\alpha_+}$ of $X_\OOO$ with $F(\alpha_-) \leq
0$ and $F(\alpha_+)\geq 0$ and $(F(\alpha_-), F(\alpha_+))\neq (0,0)$,
then a germ of $\ell$ at $z$ pointing towards $\ell(\pm\infty)$ is
contained in $\H_{\alpha_\pm}$.  
\elemm

\dem The proof of Claim 2 page 12 of \cite{BesSav12} applies.
\cqfd

\medskip
The following result implies
Theorem \ref{theo:mainintro} in the Introduction.

\btheo \label{theo:mainlourd} Let $A$ be a definite quaternion algebra over $\QQ$ and let $\OOO$ be
a maximal order in $A$. 
For every indefinite integral binary Hamiltonian form $f$ over $\OOO$,
the closest point mapping $\pi:\W_+(f)\ra\C(f)$ is a proper
$\automH$-equivariant homotopy equivalence. If the flooded
Ford-Voronoi cells for $f$ are pairwise disjoint, then the closest
point mapping $\pi:\C\W(f)\ra\C(f)$ is a $\automH$-equivariant
homeomorphism and its restriction to the waterworld $\W(f)$ is a
$\automH$-equivariant homeomorphism onto a contractible $4$-manifold
with a polyhedral boundary component homeomorphic to $\RR^3$ contained
in every flooded Ford-Voronoi cell.  
\etheo

\dem The $\automH$-equivariance properties are immediate.  We will
subdivide this proof into several steps. Unless otherwise stated,
polyhedra are compact and convex.

\medskip
\noindent {\bf Claim 1.} The closest point mapping
$\pi:\W_+(f)\ra\C(f)$ has the following properties.
\begin{enumerate}
\item The restriction of $\pi$ to any cell of $\W(f)$ is a
  homeomorphism onto its image, which is a hyperbolic polyhedron in
  the hyperbolic hyperplane $\C(f)$.
\item The restriction of $\pi$ to any flooded Ford-Voronoi cell
  $\H_\alpha$ of $f$ is a proper map onto a noncompact convex
  hyperbolic polyhedron in $\C(f)$ containing $B_\alpha(1)\cap \C(f)$
  and contained in $B_\alpha(\sqrt{D_A})\cap \C(f)$.
\item If the flooded Ford-Voronoi cells
for $f$ are pairwise disjoint, then the restriction of $\pi$ to any
cell in the boundary of a flooded Ford-Voronoi cells for $f$ is a
homeomorphism onto its image, which is a hyperbolic polyhedron in the
hyperbolic hyperplane $\C(f)$.
\end{enumerate}

\dem (1) Any $4$-dimensional cell, hence any cell, of $\W(f)$ is a hyperbolic
polyhedron in the equidistant hyperbolic hyperplane 
$$
\Scal_{\alpha,\,\beta}=\{x\in\hcr\;:\; d_\alpha(x)= d_\beta(x)\}
$$ for some $\alpha\neq \beta$ in $\PP^1_r(A)$ with $F(\alpha)F(\beta)
<0$.  Note that $\Scal_{\alpha,\,\beta}$ is not perpendicular to
$\C(f)$, otherwise $\alpha$ and $\beta$, which are the points at
infinity of a geodesic line perpendicular to $\Scal_{\alpha,\,\beta}$,
would belong to the closure of the same component of $\partial_\infty
\hcr-\C_\infty(f)$, which contradicts the fact that $F(\alpha)
F(\beta) <0$. Hence Assertion (1) of Claim 1 follows from Lemma
\ref{lem:propproj}.

\medskip\noindent(2) The closest point mapping from a horoball $H$ to a hyperbolic
hyperplane $P$ passing through the point at infinity of $H$ is a
proper map (since the intersection of $H$ with any geodesic line not
passing through its point at infinity is compact), whose image is
$H\cap P$, and which maps the geodesic segment between two points to
the geodesic segment between their images. Assertion (2) of Claim 1
hence follows from Proposition \ref{prop:propriFVcells} (2).

\medskip\noindent(3) If the flooded Ford-Voronoi cells for $f$ are pairwise disjoint,
any $4$-dimensional cell, hence any cell, in the boundary of a flooded
Ford-Voronoi cell for $f$ is a hyperbolic polyhedron in the hyperbolic
hyperplane $\Scal_{\alpha,\,\beta}$ for some $\alpha\neq \beta$ in
$\PP^1_r(A)$ with $F(\alpha)=0$ and $F(\beta)\neq 0$. Note that
$\Scal_{\alpha,\,\beta}$ is again not perpendicular to $\C(f)$,
otherwise $\alpha$ and $\beta$ would both belong to $\C_\infty(f)$,
and the Ford-Voronoi cells $\H_\alpha$ and $\H_{\beta}$ would both be
flooded for $f$ and not disjoint. The last assertion of Claim 1 follows.
\cqfd

\medskip
\noindent {\bf Claim 2.} We have the following parity properties.
\begin{enumerate}
\item Any $3$-dimensional cell $\sigma$ of $\W(f)$
not contained in a flooded Ford-Voronoi cell for $f$ belongs to an
even number of $4$-dimensional cells of $\W(f)$. 
\item
If the flooded Ford-Voronoi cells for $f$ are pairwise disjoint, then
any $3$-dimensional cell $\sigma'$ of $\W(f)$ contained in a flooded
Ford-Voronoi cell for $f$ belongs to an odd number of $4$-dimensional
cells of $\W(f)$.
\end{enumerate}

\dem (1) Since $\sigma$ has codimension $2$, the link of $\sigma$ in
the Ford-Voronoi cellulation of the manifold $\hcr$ is a
circle. Considering its intersection with the $4$-dimensional cells,
this circle subdivides into closed intervals with disjoint interiors,
each one of them contained in some Ford-Voronoi cell. By the
assumption on $\sigma$, these Ford-Voronoi cells are
nonflooded. Hence the sign of $F$ on each one of them is either $+$ or
$-$. In such a cyclic arrangement of signs, the number of sign changes
is even. Assertion (1) follows.

\medskip\noindent(2) Similarly, the link of $\sigma'$ is subdivided into at least $3$
closed intervals with disjoint interiors carrying a sign $+,0,-$.  By
the assumptions, exactly one of them, denoted by $I_0$, belongs to a
flooded Ford-Voronoi cell $\H_{\alpha_0}$ for some $\alpha_0\in
\PP^1_r(A)$, that is, carries the sign $0$. Assume for a contradiction
that the two intervals adjacent to $I_0$ carry the same sign. Let
$\beta_1,\beta_2\in\PP^1_r(A)$ be such that $\H_{\alpha_0}\cap
\H_{\beta_1}$ and $\H_{\alpha_0}\cap \H_{\beta_2}$ are the
$4$-dimensional cells corresponding to the endpoints of $I_0$. Note that
the points at $+\infty$ of the geodesic lines starting from a given
point $\alpha_0$ of $\C_\infty(f)$, passing through a geodesic line
both of whose endpoints $\beta_1,\beta_2$ are contained in the same
component $C$ of $\partial_\infty \hcr-\C_\infty (f)$ also belong to
$C$. Hence all intervals in the link of $\sigma'$ carry the same sign,
which contradicts the fact that $\sigma'$ belongs to $\W(f)$.  As for
$\sigma$, this proves that the number of sign changes between $+$ and
$-$ in the link of $\sigma'$ is odd.
\cqfd

\medskip
\noindent {\bf Claim 3.} If $\sigma$ and $\tau$ are distinct
$4$-dimensional cells of $\W(f)$ or flooded Ford-Voronoi cells for
$f$, then $\pi(\sigma)$ and $\pi(\tau)$ have disjoint interiors.

\medskip
\dem Note that no $4$-dimensional cell of $\W(f)$ is contained in a
flooded Ford-Voronoi cell for $f$.  

For a contradiction, assume that a point $p\in\C(f)$ is contained in
the interior of both $\pi(\sigma)$ and $\pi(\tau)$ and, up to moving
it a little bit, is not in the (measure $0$) image by $\pi$ of the
codimension $1$ skeleton of $X_\OOO$. Let $\ell$ be the geodesic line
through $p$ perpendicular to $\C(f)$, meeting $\sigma$ and $\tau$ at
interior points $x$ and $y$ respectively. Since the cell complex
$X_\OOO$ is locally finite, we may assume that the geodesic segment
$[x,y]$ does not meet any other $4$-dimensional cell of $\W(f)$ or
flooded Ford-Voronoi cell for $f$ than $\sigma$ and $\tau$. 

Assume for a contradiction that $[x,y]$ is contained in
$\sigma\cup\tau$. Then $\sigma$ and $\tau$ are flooded Ford-Voronoi
cells, meeting in a $4$-dimensional cell $C$, which is crossed
transversally by $[x,y]$ since $\ell$ does not meet the $3$-skeleton
of $X_\OOO$. Since $\sigma,\tau$ are flooded, their points at infinity
$\alpha,\beta\in\PP^1_r(A)$ belong to $\C_\infty(f)$. Hence the
hyperbolic hyperplane $S_{\alpha,\beta}$ equidistant to $\alpha$ and
$\beta$, which contains $\sigma$, is perpendicular to $\C(f)$. In
particular, $\ell$, which is perpendicular to $\C(f)$, is contained in
the closure of one of the two connected component of $\HH^5_\RR-
S_{\alpha,\beta}$. This contradicts the fact that $\ell$ meets
transversally $C$.

Hence $[x,y]$ is not contained in $\sigma\cup\tau$. Let $]x',y'[\;= 
[x,y]-(\sigma\cup\tau)\cap [x,y]$ with $x,x',y',y$ in this order on
$[x,y]$, so that $[x',y']$ is contained in a Ford-Voronoi cell
$\H_\alpha$ for some $\alpha\in\PP^1_r(A)$.  Let $\sigma'$ and $\tau'$
be the $4$-dimensional cells of $X_\OOO$ containing $x'$ and $y'$
respectively (note that for instance $x=x'$ and $\sigma=\sigma'$ if
$\sigma$ is a $4$-dimensional cell of $\W(f)$, but $x\neq x'$ if
$\sigma$ is a flooded Ford-Voronoi cell).

Now Lemma \ref{lem:claim3} implies that, since the two germs of the
segment $[x',y']$ at its endpoints have opposite direction, the sign
of $F(\alpha)$ should be both positive and negative, a contradiction.
\cqfd

\medskip
\noindent 
{\bf Claim 4.} The $3$-dimensional cells of the waterworld satisfy the
following properties.
\begin{enumerate}
\item No $3$-dimensional cell of $\W(f)$ is contained in two distinct
  flooded Ford-Voronoi cells.
\item Any $3$-dimensional cell $\sigma$ of $\W(f)$ not contained in a
  flooded Ford-Voronoi cell for $f$ belongs to exactly two
  $4$-dimensional cells $\tau$ and $\tau'$ of $\W(f)$, and $\pi$
  embeds their union.
\item Any $3$-dimensional cell $\sigma$ of $\W(f)$ contained in a
  flooded Ford-Voronoi cell $\H_\alpha$ for $f$ belongs to exactly one
  $4$-dimensional cell $\tau$ of $\W(f)$, and $\pi$ embeds the union
  of $\tau$ and $\tau'= \bigcup_{x\in\sigma} [x,\alpha[$.
\end{enumerate}

\dem (1) Assume for a contradiction that $\sigma$ is a $3$-dimensional
cell of $\W(f)$ contained in the flooded Ford-Voronoi cells
$\H_\alpha$ and $\H_\beta$ with $\alpha\neq \beta$ in
$\PP^1_r(A)$. Let $\tau$ be a $4$-dimensional cell of $\W(f)$
containing $\sigma$. Then the interiors of the images by $\pi$ of
$\tau$ and either $\H_\alpha$ or $\H_\beta$ are not disjoint, which
contradicts Claim 3.

Let us prove Assertions (2) and (3). Three $n$-dimensional polytopes
in $\hnr$ having a common codimension $1$ face cannot have pairwise
disjoint interiors, so that the claims on the number of
$4$-dimensional cells of $\W(f)$ containing $\sigma$ follows from
Claim 3. Since the polyhedra $\pi(\tau)$ and $\pi(\tau')$ are convex,
the result follows.  
\cqfd

\medskip
\noindent 
{\bf Claim 5.}  The $2$-dimensional cells of the waterworld satisfy the
following properties.
\begin{enumerate}
\item For every $2$-dimensional cell $\sigma$ of $\W(f)$ not contained
  in a flooded Ford-Voronoi cell for $f$, the link of $\sigma$ in
  $\W(f)$ is a circle and the union of the $4$-dimensional cells of
  $\W(f)$ containing $\sigma$ embeds in $\C(f)$ by $\pi$.
\item If the flooded Ford-Voronoi cells for $f$ are pairwise disjoint,
  for every $2$-dimensional cell $\sigma'$ of $\W(f)$ contained in a
  flooded Ford-Voronoi cell $\H_\alpha$, the link of $\sigma'$ in
  $\W(f)$ is an interval and the union of the $4$-dimensional cells of
  $\W(f)$ containing $\sigma'$ and of the geodesic rays $[x,\alpha[$
  for $x$ in the two $3$-dimensional cells of $\W(f)\cap \partial
  \H_\alpha$ containing $\sigma'$ embeds in $\C(f)$ by $\pi$.
\end{enumerate}

\dem (1) By Claim 4, the link $Lk(\sigma)$ of $\sigma$ in $\W(f)$ is a
disjoint union of circles. Each component of $Lk(\sigma)$ corresponds
to a finite set of $4$-dimensional cells cyclically arranged around
$\sigma$. By Claim 4 again, their images by $\pi$ are not folded,
hence are cyclically arranged around $\pi(\sigma)$. If $Lk(\sigma)$
was not connected, the image of two $4$-dimensional cells of $\W(f)$
by $\pi$ would have intersecting interiors, contradicting Claim 3. 

\medskip\noindent(2) An analogous proof gives that the link of $\sigma'$ in $\C\W(f)$
is a circle.  \cqfd

\medskip
\noindent 
{\bf Claim 6.} The $1$-dimensional cells of the waterworld satisfy the
following properties.
\begin{enumerate}
\item For every $1$-dimensional cell $\sigma$ of $\W(f)$ not
contained in a flooded Ford-Voronoi cell for $f$, the link of $\sigma$
in $\W(f)$ is a $2$-sphere and the union of the $4$-dimensional cells
of $\W(f)$ containing $\sigma$ embeds in $\C(f)$ by $\pi$.
\item If the flooded Ford-Voronoi cells for $f$ are pairwise disjoint,
  for every $1$-dimensional cell $\sigma'$ of $\W(f)$ contained in a
  flooded Ford-Voronoi cell $\H_\alpha$, the link of $\sigma'$ in
  $\W(f)$ is a $2$-disc and the union of the $4$-dimensional cells of
  $\W(f)$ containing $\sigma'$ and of the geodesic rays $[x,\alpha[$
      for $x$ in any $3$-cell of $\W(f)\cap\partial \H_\alpha$
      containing $\sigma'$ embeds in $\C(f)$ by $\pi$.
\end{enumerate}

\dem (1) By Claim 5, the links of the vertices of the link
$Lk(\sigma)$ of $\sigma$ in $\W(f)$ are circles, hence $Lk(\sigma)$ is
a compact surface, mapping locally homeomorphically to
$Lk(\pi(\sigma))$ by $\pi$, which is a $2$-sphere. Hence
$Lk(\pi(\sigma))$ is a union of $2$-spheres, again with only one of
them by Claim 3.

\medskip\noindent(2) The proof that the link of $\sigma'$ in $\C\W(f)$ is a $2$-sphere
is similar. \cqfd

\medskip
\noindent 
{\bf Claim 7.} The vertices of the waterworld satisfy the
following properties. 
\begin{enumerate}
\item For every vertex $v$ of $\W(f)$ not contained in a flooded
  Ford-Voronoi cell for $f$, the link of $v$ in $\W(f)$ is a
  $3$-sphere and the union of the $4$-dimensional cells of $\W(f)$
  containing $v$ embeds in $\C(f)$ by $\pi$.
\item If the flooded Ford-Voronoi cells for $f$ are pairwise disjoint,
  for every vertex $v'$ of $\W(f)$ contained in a flooded Ford-Voronoi
  cell $\H_\alpha$, the link of $v'$ in $\W(f)$ is a $3$-disc and the
  union of the $4$-dimensional cells of $\W(f)$ containing $v'$ and of
  the geodesic rays $[x,\alpha[$ for $x$ in any $3$-cell of
      $\W(f)\cap\H_\alpha$ containing $v'$ embeds in $\C(f)$ by $\pi$.
\end{enumerate}

\dem The proof is similar to the previous one.
\cqfd

\medskip
Now, the properness of $\pi:\W_+(f)\ra \C(f)$ follows from the fact
that $\pi$ is $\automH$-equivariant, that $\automH$ acts cocompactly
on $\W(f)$ and with finitely many orbits on the set of flooded
Ford-Voronoi cells by Proposition
\ref{prop:waterworldmodautomcompact}, and from its properness when
restricted to each flooded Ford-Voronoi cell (see Claim 1).

Claim 7 proves that when the flooded Ford-Voronoi cells for $f$ are
pairwise disjoint, the map $\pi:\C\W(f)\ra \C(f)$ is a proper local
homeomorphism betwen locally compact spaces, hence is a covering
map. Since $\C(f)$ is simply connected, $\pi$ is hence a homeomorphism
on each of the connected components of $\C\W(f)$. But since $\pi$ is
injective outside the codimension $1$ skeleton by Claim 3, it follows
that $\C\W(f)$ is connected and $\pi$ is a homeomorphism. This
concludes the proof of Theorem \ref{theo:mainlourd}.
\cqfd

\appendix
\section{An algebraic description of the distance to the cusps}
\label{sec:appAalgebraicdistcusp}

%%%%%%%%%%%%%%%%%%%
%% Modified 30.10.19 by Frédéric
%%%%%%%%%%%%%%%%%%%
Let $A$ be a definite quaternion algebra over $\QQ$ and let $\OOO$ be
a maximal order in $A$.  In this independent appendix, following
Mendoza \cite{Mendoza80} in the Hermitian case, we give an algebraic
description of the distance functions $d_\alpha$ to the rational
points at infinity $\alpha\in\PP^1_r(A)$, defined just before
Proposition \ref{prop:propridalpha}.

An {\it $\OOO$-flag} is a right $\OOO$-submodule $L$ of the right
$\OOO$-module $\OOO\times\OOO$, with rank one (that is, $LA$ is a line
in the $A$-vector space $A\times A$), such that the quotient
$(\OOO\times \OOO)/L$ has no torsion. We denote by $\F_\OOO$ the set
of $\OOO$-flags.

For all right $\OOO$-submodules $M$ of $A\times A$ and $v\in
A\times A-\{0\}$, let us define
$$
M_v=\{x\in A\;:\; vx\in M\}\;.
$$
Note that for every $\lambda\in A-\{0\}$, we immediately have
\begin{equation}\label{eq:homothetyMv}
\lambda M_{v\lambda}=M_v\;.
\end{equation}

\bexem 
Recall that the {\it inverse} $I^{-1}$ of a left fractional
ideal $I$ of $\OOO$ is the right fractional ideal of $\OOO$
$$
I^{-1}=\{x\in A\;:\; IxI\subset I\}\;.
$$
It is well known and easy to check that for every $a,b\in\OOO$, if
$ab\neq 0$, then
\begin{equation}\label{eq:invsomeginter}
(\OOO a+\OOO b)^{-1}=a^{-1}\OOO \cap b^{-1}\OOO \;.
\end{equation}
We claim that if $v=(a,b)$, then
\begin{equation}\label{eq:calculMv}
(\OOO\times\OOO)_v=(\OOO a+\OOO b)^{-1}\;.
\end{equation}
Indeed, if $a,b\neq 0$, then by Equation \eqref{eq:invsomeginter}
$$
(\OOO\times\OOO)_v=\{x\in A\;:\; (ax,bx)\in \OOO\times \OOO\}=
a^{-1}\OOO \cap b^{-1}\OOO=(\OOO a+\OOO b)^{-1}\;.
$$
The result is immediate if $a=0$ or $b=0$.
\eexem 

\medskip

\bprop (1) For all right $\OOO$-submodule $M$ of $A\times A$ and $v\in
A\times A-\{0\}$, the subset $M_v$ of $A$ is a right fractional ideal
of $\OOO$.

\smallskip\noindent (2) For every $v\in A\times A-\{0\}$, the subset
$v(\OOO\times\OOO)_{v}$ of $\OOO\times \OOO$ is an $\OOO$-flag.

\smallskip\noindent (3) For all $\OOO$-flags $L$ and all $v\in
L-\{0\}$, we have
$$
L=v(\OOO\times\OOO)_v\;.
$$

\smallskip\noindent (4) The map $\SL_2(A)\times \F_\OOO\ra \F_\OOO$ defined by
$$
(g,L)\mapsto (gv)(\OOO\times\OOO)_{gv}
$$
for any $v\in L-\{0\}$ is an action on the set $\F_\OOO$ of
$\OOO$-flags of the group $\SL_2(A)$.

\smallskip\noindent (5) The map $\Theta': \PP^1_r(A)\ra \F_\OOO$
defined by $[a:b]\mapsto (a,b)(\OOO\times\OOO)_{(a,b)}$ is a
$\SL_2(A)$-equivariant bijection. 
\eprop

\dem 
(1) This follows immediately from the fact that $M$ is stable by
addition and by multiplications on the right by the elements of
$\OOO$.

\medskip\noindent 
(2) Let $L=v(\OOO\times\OOO)_{v}\subset vA$. Then $L$ is contained in
$\OOO\times \OOO$ by the definition of $(\OOO\times \OOO)_{v}$ and is
a right $\OOO$-submodule of $\OOO\times \OOO$ by Assertion (1). Since
$v\neq 0$, note that $(\OOO\times\OOO)_{v}$ is a nonzero right
fractional ideal, so that $L\neq 0$ and $L$ has rank one.

Assume that $w\in\OOO\times\OOO$ has its image in $(\OOO\times\OOO)/L$
which is torsion. Then there exists $y\in\OOO-\{0\}$ and $x\in A$ such
that $wy=vx$. Hence $w=vxy^{-1}$. Since $w\in\OOO\times\OOO$, this
implies that $xy^{-1}\in (\OOO\times\OOO)_{v}$, so that $w\in L$, and
the image of $w$ in $(\OOO\times\OOO)/L$ is zero.

\medskip\noindent 
(3) As $L$ has rank $1$ and $v\in L-\{0\}$, we have $L\subset vA \cap
(\OOO\times\OOO)=v(\OOO\times\OOO)_v$.

Conversely, for every $x\in (\OOO\times\OOO)_v$ so that $vx\in \OOO
\times \OOO$, let us prove that $vx\in L$. Since $x\in A$ which is the
field of fractions of $\OOO$, there exists $y\in \OOO$ such that
$xy\in\OOO$. Hence $(vx)y=v(xy)$ belongs to $L$, since $v\in L$ and
$L$ is a right $\OOO$-module. In particular, the image of $vx$ in
$(\OOO\times\OOO)/L$ is torsion. Since $L$ is an $\OOO$-flag, this
implies that this image is zero, as wanted.  This proves that
$v(\OOO\times\OOO)_{v}$ is contained in $L$, hence is equal to $L$ by
the previous inclusion.

\medskip\noindent
(4) Let us prove that this map is well defined. If $v,w\in L-\{0\}$,
since $L$ has rank one, there exists $x\in A-\{0\}$ such that $w=vx$,
thus, for every $g\in \SL_2(A)$, by the linearity on the right of
$g$ and by Equation \eqref{eq:homothetyMv}, we have
$$
(gw)(\OOO\times\OOO)_{gw}=(gv)x(\OOO\times\OOO)_{(gv)x}=
(gv)(\OOO\times\OOO)_{gv}\;.
$$
The fact that this map is an action is then immediate: for all
$g,g'\in \SL_2(A)$ and $L\in \F_\OOO$, let $v\in L-\{0\}$ and
$\lambda\in A$ be such that $gv\lambda\in (gv)(\OOO\times\OOO)_{gv}
-\{0\}$; then using twice Equation \eqref{eq:homothetyMv} and the linearity,
we have
\begin{align*}
g'(gL)&=g'\big(gv(\OOO\times\OOO)_{gv}\big)
=g'\big(gv\lambda(\OOO\times\OOO)_{gv\lambda}\big)
=g'(gv\lambda)(\OOO\times\OOO)_{g'(gv\lambda)}\\ & =
(g'g)v\lambda(\OOO\times\OOO)_{(g'g)v\lambda}=
(g'g)v(\OOO\times\OOO)_{v}= (g'g)L\;.
\end{align*}

\medskip\noindent
(5) For every $\alpha=[a:b]\in\PP^1_r(A)$, the subset $(a,b)
(\OOO\times\OOO)_{(a,b)}$, which is an $\OOO$-flag by Assertion (2),
does not depend on the choice of homogeneous coordinates of $\alpha$
by Equation \eqref{eq:homothetyMv}. Hence the map $\Theta'$ is well
defined, and equivariant by the definition of the action of $\SL_2(A)$
on $\F_\OOO$.

The fact that $\Theta'$ is onto follows from Assertion (3). 
Clearly, it is  one-to-one since if $(a,b) (\OOO\times\OOO)_{(a,b)}=(c,d)
(\OOO\times\OOO)_{(c,d)}$, then there is $\lambda\in A-\{0\}$ such
that $(a,b)=(c,d)\lambda$.
\cqfd

\medskip 
Let $f:\HH\times\HH\ra\RR$ be a positive definite binary Hamiltonian
form and let $L$ be a rank one right $\OOO$-submodule of $\OOO\times \OOO$. Then
$L$ is a rank $4$ free $\ZZ$-submodule of $\HH\times \HH$, and we
denote by $\langle L\rangle_\RR$ the $4$-dimensional real vector
subspace of $\HH\times \HH$ generated by $L$, endowed with the
restriction of the scalar product $\langle\cdot,\cdot\rangle_f$ on
$\HH\times \HH$ defined by $f$, hence with the induced volume
form. Recall that for all $z,z'\in\HH\times\HH$, we have
\begin{equation}\label{eq:scalarproductf}
\langle z,z'\rangle_f =
\frac{1}{2}\big( f(z+z')-f(z)-f(z')\big)\;.
\end{equation}

We define the covolume of $L$ for $f$ as
$$
\covol_f L =\Vol(\langle L\rangle_\RR/L)\;.
$$
Recall
that if $G=(\langle e_i,e_j\rangle_f)_{1\leq i,j\leq 4}$ is the Gram
matrix of a $\ZZ$-basis $(e_1,e_2,e_3,e_4)$ of $L$ for the scalar
product $\langle\cdot,\cdot\rangle_f$, then
\begin{equation}\label{eq:gram}
\covol_f L =(\det G)^{\frac{1}{2}}\;.
\end{equation}
See for instance \cite[Vol~2, prop.~8.11.6]{Berger79}.

\btheo \label{theo:computdistcuspalg}
For all $x\in\hcr$ and $\alpha\in \PP^1_r(A)$, we have
$$
d_\alpha(x)=\frac{2}{\sqrt{D_A}}\;
\big(\covol_{\Theta(x)}\Theta'(\alpha)\big)^{\frac{1}{2}}\;.
$$
\etheo

\dem 
Fix $a,b\in\OOO$ such that $\alpha=[a:b]$. Let $f=\Theta(x)$, $L=
\Theta'(\alpha)=(a,b)(\OOO\times\OOO)_{(a,b)}$ and $L'= (a,b)\OOO$.
Since $a,b\in\OOO$, we have $\OOO\subset (\OOO\times\OOO)_{(a,b)}$,
hence $L'$ is a finite index $\ZZ$-submodule in $L$. Furthermore, by
Equation \eqref{eq:calculMv} and the relation (see Equation
\eqref{eq:rednormindex}) between the norm and reduced norm of a left
integral ideal of $\OOO$, we have
\begin{align}
[L:L']&=[(\OOO\times\OOO)_{(a,b)}:\OOO]=[(\OOO a+\OOO b)^{-1}:\OOO]=
[\OOO:\OOO a+\OOO b]\nonumber\\ &
 \label{eq:computindice} =\n(\OOO a +\OOO b)^2\;.
\end{align}

Let $(x_1,x_2,x_3,x_4)$ be a $\ZZ$-basis of $\OOO$, so that
$(\,(a,b)x_i)_{1\leq i\leq 4}$ is a $\ZZ$-basis of $L'$.  Using
Equation \eqref{eq:scalarproductf} and the fact that $f((u,v)\lambda)
= \n(\lambda)f(u,v)$ for all $u,v,\lambda\in\HH$, we have for $1\leq
i,j\leq 4$,
\begin{align*}
\langle(a,b)x_i,(a,b)x_j\rangle_f & =\frac{1}{2}
\big( f\big((a,b)(x_i+x_j)\big)-f((a,b)\,x_i)-f((a,b)\,x_j)\big)\\
& = \frac{f(a,b)}{2}\,(\n(x_i+x_j)-\n(x_i)-\n(x_j))=
\frac{f(a,b)}{2}\,\tr(\,\overline{x_i}\,x_j)\;.
\end{align*}
Note that $(u,v)\mapsto \frac{1}{2}\tr(\,\overline{u}\,v)$ is the
standard Euclidean scalar product on $\HH$ (making the standard basis
$(1,i,j,k)$ orthonormal), hence $\big(\frac{1}{2} \tr(\,\overline{x_i} 
\,x_j)\big)_{1\leq i,j\leq 4}$ is the Gram matrix of the $\ZZ$-lattice 
$\OOO$ in the Euclidean space $\HH$. Therefore, by Equation 
\eqref{eq:gram} and by \cite[Lem.~5.5]{KraOse90}, we have
\begin{equation}\label{eq:gramclassique}
\big(\det\big(\tr(\,\overline{x_i}\, x_j)\big)_{1\leq i,j\leq 4} 
\big)^{\frac{1}{2}}= (2^4)^{\frac{1}{2}}\;\Vol(\HH/\OOO)= 4 \frac{D_A}{4}
=D_A\;.
\end{equation}
Thus using Equations \eqref{eq:gram}, \eqref{eq:computindice}
and \eqref{eq:gramclassique}, we have
\begin{align*}
\covol_f(L)&=\frac{1}{[L:L']}\;\covol_f(L')=
\frac{1}{[L:L']}\big(\det\big(\langle(a,b)x_i,
(a,b)x_j\rangle_f\big)_{1\leq i,j\leq 4}\big)^{\frac{1}{2}}\\ 
& =\frac{1}{[L:L']}\big(\frac{f(a,b)}{2}\big)^2
\big(\det\big(\tr(\,\overline{x_i}\, x_j)\big)_{1\leq i,j\leq 4} 
\big)^{\frac{1}{2}}
= \frac{D_A}4 \;\frac{f(a,b)^2}{\n(\OOO a +\OOO b)^2}\;.
\end{align*}
By Proposition \ref{prop:propridalpha} (2), this proves
Theorem \ref{theo:computdistcuspalg}.
\cqfd

{\small \bibliography{../biblio} }
%{\small \bibliography{/users/jouniparkkonen/Dropbox/Latex//viitteet} }

\bigskip
{\small
\noindent \begin{tabular}{l} 
Department of Mathematics and Statistics, P.O. Box 35\\ 
40014 University of Jyv\"askyl\"a, FINLAND.\\
{\it e-mail: jouni.t.parkkonen@jyu.fi}
\end{tabular}
\medskip

\noindent \begin{tabular}{l}
Laboratoire de math\'ematique d'Orsay,\\
UMR 8628 Univ. Paris-Sud et CNRS,\\
Universit\'e Paris-Saclay,\\
91405 ORSAY Cedex, FRANCE\\
{\it e-mail: frederic.paulin@math.u-psud.fr}
\end{tabular}
}

\end{document}